\documentstyle{amltd}
\begin{document}
\annalsline{157}{2003}
\received{February 20, 2001}
\startingpage{807}
\def\bye{\end{document}}
 \font\tenrm=cmr10
\def\ritem#1{\item[{\rm #1}]}
\def\pref#1{({\ref{#1}})}
\catcode`\@=11
\font\twelvemsb=msbm10 scaled 1100
\font\tenmsb=msbm10
\font\ninemsb=msbm10 scaled 800
\newfam\msbfam
\textfont\msbfam=\twelvemsb  \scriptfont\msbfam=\ninemsb
  \scriptscriptfont\msbfam=\ninemsb
\def\msb@{\hexnumber@\msbfam}
\def\Bbb{\relax\ifmmode\let\next\Bbb@\else
 \def\next{\errmessage{Use \string\Bbb\space only in math
mode}}\fi\next}
\def\Bbb@#1{{\Bbb@@{#1}}}
\def\Bbb@@#1{\fam\msbfam#1}
\catcode`\@=12

 \catcode`\@=11
\font\twelveeuf=eufm10 scaled 1100
\font\teneuf=eufm10
\font\nineeuf=eufm7 scaled 1100
\newfam\euffam
\textfont\euffam=\twelveeuf  \scriptfont\euffam=\teneuf
  \scriptscriptfont\euffam=\nineeuf
\def\euf@{\hexnumber@\euffam}
\def\frak{\relax\ifmmode\let\next\frak@\else
 \def\next{\errmessage{Use \string\frak\space only in math
mode}}\fi\next}
\def\frak@#1{{\frak@@{#1}}}
\def\frak@@#1{\fam\euffam#1}
\catcode`\@=12

\def\a{\alpha}
\def\b{\beta}
\def\d{\delta}
\def\e{\varepsilon}
\def\ve{\varepsilon}
\def\f{\phi}
\def\vf{\varphi}
\def\g{\gamma}
\def\k{\kappa}
\def\l{\lambda}
\def\r{\rho}
\def\s{\sigma}
\def\t{\tau}
\def\th{\theta}
\def\vp{\varphi}
\def\z{\zeta}
\def\o{\omega}
\def\D{\Delta}
\def\L{\Lambda}
\def\G{\Gamma}
\def\O{\Omega}
\def\S{\Sigma}
\def\Th{\Theta}
\def\del #1{\frac{\partial^{#1}}{\partial\l^{#1}}}

\def\h{\eta}

\def\E{{\Bbb E}}
\def\N{{\Bbb  N}}
\def\M{{\Bbb M}}
\def\R{{\Bbb R}}
\def\Z{{\Bbb Z}}
\def\1{{1\kern-.25em\hbox{\rm I}}}
\def\eu{{1\kern-.25em\hbox{\sm I}}}

\def\C{{\Bbb C}}
\def\P{{\Bbb P}}
\def\eop{{ \vrule height7pt width7pt depth0pt}\par\bigskip}

\def\dert{\frac{\partial}{\partial t}}
\def\der{\frac{d}{dx}}
\def\del{\partial}
\def\tr{\hbox{tr}}
\def\pb{\Phi_\b}
\def\AA{{\cal A}}
\def\BB{{\cal B}}
\def\CC{{\cal C}}
\def\DD{{\cal D}}
\def\EE{{\cal E}}
\def\FF{{\cal F}}
\def\GG{{\cal G}}
\def\II{{\cal I}}
\def\JJ{{\cal J}}
\def\KK{{\cal K}}
\def\LL{{\cal L}}
\def\OO{{\cal O}}
\def\PP{{\cal P}}
\def\SS{{\cal S}}
\def\TT{{\cal T}}
\def\NN{{\cal N}}
\def\MM{{\cal M}}
\def\WW{{\cal W}}
\def\VV{{\cal V}}
\def\UU{{\cal U}}
\def\LL{{\cal L}}
\def\XX{{\cal X}}
\def\ZZ{{\cal Z}}
\def\RR{{\cal R}}
\def\QQ{{\cal Q}}
\def\A{{\cal A}}
\def\T{\RR}
\def\chap #1#2{\line{\ch #1\hfill}\numsec=#2\numfor=1}
\def\sign{\,\hbox{sign}\,}
\def\un #1{\underline{#1}}
\def\ov #1{\overline{#1}}
\def\ovm{{\overline m}}
\def\ba{{\backslash}}
\def\sb{{\subset}}
\def\wt{\widetilde}
\def\wh{\widehat}
\def\sp{{\supset}}
\def\rar{{\rightarrow}}
\def\lar{{\leftarrow}}
\def\em{{\emptyset}}
\def\inn{\hbox{int}\,}
\def\sminn{\hbox{\ftn int}\,}
\def\smcl{\hbox{\ftn cl}\,}
\def\rr#1{r^*(M,#1)}

\def\Pp#1{\Psi_{\b,N,\d}(#1)}
\def\pp#1{\phi_{\b,\d}(#1)}
\def\adb{a_{\d,\b}}
\def\aa{a}
\def\ee{\frac{\xi^t\xi}{N}}
\def\Rr#1{r(M,#1)}
\def\rr#1{r^*(M ,#1)}
\def\pat#1{\xi^{#1}}
\def\ov#1{\overline{#1}}
\def\jg{J_\g}
\def\sqr#1#2{{\vcenter{\vbox{\hrule height.#2pt
     \hbox{\vrule width.#2pt height#1pt \kern#1pt
   \vrule width.#2pt}\hrule height.#2pt}}}}
\def\qed{ $\mathchoice\sqr64\sqr64\sqr{2.1}3\sqr{1.5}3$} 

\def\sfrac#1#2{{\textstyle{#1\over #2}}}
\def\text#1{\quad{\hbox{#1}}\quad}
\font\pr=cmbxsl10
\font\thbf=cmbxsl10 scaled\magstephalf

\def\ds{\displaystyle}\def\st{\scriptstyle}\def\sst{\scriptscriptstyle}

\def\ni{\noindent}

\title{Constrained steepest descent\\ in the 2-Wasserstein metric} 
\shorttitle{The 2-Wasserstein metric} 

  \acknowledgements{The work of the first named author was partially supported by
U.S.\ N.S.F.\ grant DMS-00-70589.  The work of the second named author was partially supported by
U.S.\ N.S.F.\ grants DMS-99-70520 and DMS-00-74037.}
 \twoauthors{E.\ A.\ Carlen}{W.\ Gangbo}
 \institutions{Georgia Institute of Technology,
  Atlanta, GA\\
{\eightpoint {\it E-mail addresses\/}: carlen@math.gatech.edu}\\ 
 \hglue.97in  {\eightpoint  gangbo@math.gatech.edu}
\vglue-12pt}
\vfil
\centerline{\bf Abstract } 
\vglue12pt

We study several constrained variational problems in the $2$-Wasserstein metric
for which the set of probability densities satisfying the constraint is not 
closed. 
For example, given a probability 
density $F_0$ on $\R^d$ and a time-step\break $h>0$, we seek to minimize
$I(F) = hS(F) + W_2^2(F_0,F)$
over all of the probability densities $F$ that have the same mean and variance 
as $F_0$, where $S(F)$ is the entropy of $F$. We prove existence of minimizers.
We also analyze the induced geometry of the set of densities satisfying the 
constraint on the variance and means, and we determine all of the geodesics on 
it. From this, we determine a criterion for convexity of functionals
in the induced geometry. It turns out, for example, that the entropy is
uniformly strictly convex on the constrained manifold, though not uniformly 
convex without the constraint.
The problems solved here arose in 
a study of a variational approach to 
constructing and studying solutions of the nonlinear kinetic Fokker-Planck 
equation, which is 
briefly described here
and fully developed in a companion paper. 
 \vfil

\def\sni#1{\smallbreak\noindent{#1}. }
\centerline{\bf Contents}
\sni{1} Introduction
\sni{2} Riemannian geometry of the $2$-Wasserstein metric
\sni{3} Geometry of the constraint manifold
\sni{4} The Euler-Lagrange equation
\sni{5} Existence of minimizers
\smallbreak\noindent \phantom{1. } References
\pagebreak

\section{Introduction}

Recently there has been considerable progress in understanding a wide range of 
dissipative evolution equations  
in terms of variational problems involving the Wasserstein metric.
In particular, Jordan, Kinderlehrer and Otto, have shown in \cite{JKO}
that the heat equation
is gradient flow for the entropy functional in the $2$-Wasserstein metric. We can 
arrive most rapidly 
to the point of departure for our own problem, which concerns constrained 
gradient flow, by reviewing this result.

Let ${\cal P}$ denote the set of probability densities  on $\R^d$ with finite 
second moments; i.e.,
the set of all nonnegative measurable functions $F$ on $\R^d$ such that
$\int_{\R^d}F(v){\rm d}v = 1$ and $\int_{\R^d}|v|^2F(v){\rm d}v <\infty$. We use 
$v$ and $w$ to denote points in $\R^d$
since in the problem to be described below they represent velocities. 
Equip ${\cal P}$  with the $2$-Wasserstein metric, $W_2(F_0,F_1)$, 
where
\begin{equation}
W_2^2(F_0,F_1) = \inf_{\gamma\in{\cal C}(F_0,F_1)}
\int_{\R^d\times\R^d}
{1\over 2}|v-w|^2\gamma({\rm d}v,{\rm d}w)\ .\label{infimum}
\end{equation}
Here, ${\cal C}(F_0,F_1)$ consists of all {\it couplings} of $F_0$ and 
$F_1$; i.e.,
all probability measures 
$\gamma$ on $\R^d\times\R^d$ such that for all test functions $\eta$ on $\R^d$
$$\int_{\R^d\times\R^d}\eta(v)\gamma({\rm d}v,{\rm d}w) = 
\int_{\R^d}\eta(v)F_0(v){\rm d}v$$
and
$$\int_{\R^d\times\R^d}\eta(w)\gamma({\rm d}v,{\rm d}w) = 
\int_{\R^d}\eta(w)F_1(w){\rm d}v\ .$$

The infimum
in \pref{infimum} is actually a minimum, and it is attained at a unique
point $\gamma_{F_0,F_1}$ in ${\cal C}(F_0,F_1)$. Brenier \cite{Bren} was able 
to characterize
this unique minimizer, and then further results of
Caffarelli \cite{Caff92}, Gangbo \cite{Gangbo} and  McCann 
\cite{McCann95} shed considerable light 
on the nature of this minimizer.

Next, let the entropy $S(F)$ be defined by
\begin{equation}
S(F) = \int_{\R^d}F(v)\ln F(v){\rm d}v\ .\label{xent} \end{equation}
This is well defined, with $\infty$ as a possible value, 
since $\int_{\R^d}|v|^2F(v){\rm d}v<\infty$.
\vglue4pt 
The following scheme for 
solving the linear heat equation was  introduced in \cite{JKO}: Fix an 
initial density $F_0$ with 
$\int_{\R^d}|v|^2F_0(v){\rm d}v$ finite, and also fix a time step $h>0$.
Then inductively define $F_k$ in terms of $F_{k-1}$ by choosing $F_k$ to 
minimize the functional
\begin{equation}
F \rightarrow \left[W_2^2(F_{k-1},F) + hS(F)\right]
\label{heatfunc}
\end{equation}
on ${\cal P}$.
It is shown in \cite{JKO} that there is a unique minimizer $F_k\in {\cal P}$, 
so that each $F_k$ is well 
defined. Then the time-dependent probability density $F^{(h)}(v,t)$
is defined by putting 
$ F^{(h)}(v,kh) = F_k$
and interpolating when $t$ is not an integral multiple of $h$. 
Finally, it is shown that for each $t$
$F(\cdot,t) = \lim_{h\to 0}F^{(h)}(\cdot,t)$
exists weakly in $L^1$, and that the resulting time-dependent probability
density solves the heat equation
$\partial/\partial tF(v,t) = 
\Delta F(v,t)$
with 
$\lim_{t\to 0}F(\cdot, t) = F_0$. 

This variational approach is particularly useful when the functional being 
minimized with each time step
is convex in the geometry associated to the\break 2-Wasserstein metric. It makes sense 
to 
speak of convexity in
this context
since, as McCann showed \cite{McCann95}, when ${\cal P}$  is equipped with 
the $2$-Wasserstein metric, 
 every pair of elements $F_0$ and $F_1$ is connected by a unique continuous path 
$t\mapsto F_t$,
$0\le t \le 1$, such that
$W_2(F_0,F_t) + W_2(F_t,F_1) = W_2(F_0,F_1)$
for all such~$t$. It is natural to refer to this path as the geodesic connecting 
$F_0$ and $F_1$, and
we shall do so. A functional $\Phi$ on ${\cal P}$ is {\it displacement convex} 
in McCann's sense if
$t\mapsto\Phi(F_t)$ is convex on $[0,1]$
for every $F_0$ and $F_1$ in ${\cal P}$. It turns out that the entropy $S(F)$ is 
a convex function of $F$ in this sense.

Gradient flows of convex functions in Euclidean space are well known to have 
strong contractive properties, 
and Otto \cite{Ottop}
showed that the same is true in ${\cal P}$, and applied this to obtain strong 
new results on 
rate of relaxation of certain solutions of the porous medium equation.

Our aim is to extend this line of analysis to a range of problems that are not 
purely dissipative, but 
which also satisfy certain {\it conservation laws}. An important example of such 
an evolution 
is given by the Boltzmann equation  
$$
{\partial\over \partial t}f(x,v,t) + \nabla_x\cdot \left(vf(x,v,t)\right) = 
{\cal Q}\left(f\right)(x,v,t)$$  
where for each $t$, $f(\cdot,\cdot,t)$ is a probability density on the phase 
space ${\Lambda\times \R^d}$
of a molecule in a region $\Lambda\subset \R^d$, and ${\cal Q}$ is a nonlinear 
operator representing the effects of 
collisions to the evolution of molecular velocities. This evolution is 
dissipative and decreases the entropy
while formally conserving the energy $\int_{\Lambda\times \R^d}|v|^2f(x,v,t){\rm 
d}x{\rm d}v$ and the
momentum $\int_{\Lambda\times \R^d}vf(x,v,t){\rm d}x{\rm d}v$. A good deal is 
known about this equation
\cite{Cer}, but there is not yet an existence theorem for solutions that 
conserve the energy, 
nor is there any general uniqueness result. 

The investigation in this paper arose in the
study of a related equation, the nonlinear kinetic Fokker-Planck equation 
to which we have applied an analog of the scheme in \cite{JKO} to the 
evolution of the
conditional probability densities  
$F(v;x)$ for the velocities of the molecules at $x$; i.e., for the contributions 
of the 
collisions to the evolution 
of the distribution of velocities of particles in a gas. These collisions are 
supposed to conserve both
the ``bulk velocity'' $u$ and ``temperature'' $\theta$, of the distribution 
where
\begin{equation}
u(F) = \int_{\R^d}vF(v){\rm d}v\qquad{\rm and}\qquad  \theta(F) = {1\over 
d}\int_{\R^d}|v|^2F(v){\rm d}v. \label{fj4} \end{equation}
For this reason we add a constraint to the variational problem in 
\cite{JKO}. 
Let $u \in \R^d$ and $\theta > 0$ be given.
Define the subset 
${\cal E}_{u,\theta}$ of ${\cal P}$
specified by
\begin{equation}
{\cal E}_{u,\theta} = \left\{F\in{\cal P}\ \biggl|\ {1\over 
d}\int_{\R^d}|v-u|^2 
F(v){\rm d}v = \theta
\quad{\rm and}\quad \int_{\R^d}vF(v){\rm d}v = u\ \right\}\ .\label{mmdef} \end{equation}
This is the set of all probability
densities with a mean $u$ and a variance $d\theta$, and we use ${\cal E}$ to 
denote it 
because the constraint on the variance is interpreted as an internal energy 
constraint in the context discussed above.

Then given $F_0\in {\cal E}_{u,\theta}$, define the functional $I(F)$ on ${\cal 
E}_{u,\theta}$ by
\begin{equation}I(F) = \left[{W_2^2(  F_0,F)
\over \theta} +hS(F)\right]\ .\label{ffudef} \end{equation}
Our main goal is to study the minimization problem associated with determining
\begin{equation}\inf\left\{ I(F)\ \big|\ F\in{\cal E}_{u,\theta}\right\}\ .\label{zmain} \end{equation}
Note that this problem is scale invariant in that if $F_0$ is rescaled, the 
minimizer $F$
will be rescaled in the same way, and in any case, this normalization, with 
$\theta$ in the denominator,
 is dimensionally natural.

Since the constraint is {\it not weakly closed}, existence of minimizers does 
not follow as easily as in the 
unconstrained case. The same difficulty arises in the determination of the 
geodesics in ${\cal E}_{u,\theta}$.

We build on previous work on the geometry of ${\cal P}$ in the $2$-Wasserstein 
metric, and Section 2 contains
a brief exposition of the relevant results. While this section is largely 
review, several of the simple proofs given here 
do not seem to be in the literature, and are more readily adapted to the 
constrained setting.

In Section 3, we analyze the geometry of ${\cal E}$, and determine its 
geodesics. 
As mentioned above, since ${\cal E}$ is not weakly closed, direct methods do not 
yield the geodesics.
The characterization of the geodesics is quite explicit, and from it we deduce
a criterion for convexity in ${\cal E}$, and show that the entropy is uniformly 
strictly convex, in contrast with the
unconstrained case. 

In Section 4, we turn to the variational problem \pref{zmain}, and
determine the Euler-Lagrange equation associated with it, and several 
consequences of the Euler-Lagrange equation.

In Section 5 we introduce  a  variational problem that is dual to \pref{zmain}, and 
by analyzing it, we produce a  minimizer for $I(F)$. We conclude the paper in 
Section 6 by discussing some open problems and possible applications.

We would like to thank Robert McCann and Cedric Villani for many enlightening discussions
on the subject of mass transport. We would also like to thank the referee, 
whose questions and suggestions have lead us 
to clarify   the exposition significantly.

\section{Riemannian geometry of the $2$-Wasserstein metric}

The purpose of this section is to collect a number of facts concerning the\break
$2$-Wasserstein metric and its associated Riemannian geometry. The Riemannian 
point of view has been developed by several authors, prominently including 
McCann, Otto, and Villani. Though for the most part the facts presented in this 
section are known, there is no single convenient reference for all of them.  
Moreover, it seems that some 
of the proofs 
and formulae that we use do not appear elsewhere in the literature.

We begin by recalling the identification of the geodesics in ${\cal P}$ equipped 
with the $2$-Wasserstein metric. The fundamental facts from which we start are 
these: The infimum in \pref{infimum} is actually a minimum, and it is attained at 
a unique point $\gamma_{F_0,F_1}$ in ${\cal C}(F_0,F_1)$, and this 
measure is such that there exists a pair of dual convex functions
$\phi$ and $\psi$ such that for all bounded measurable functions $\eta$ 
on $\R^d\times\R^d$,
\begin{eqnarray} \int_{\R^d\times\R^d}\eta(v,w)\gamma_{F_0,F_1}({\rm d}v,{\rm d}w)
&= &\int_{\R^d}\eta(v,\nabla\phi(v))F_0{\rm d}v\label{apform}\\
& =& 
\int_{\R^d}\eta(\nabla\psi(w),w)F_1{\rm d}w\ . \nonumber\end{eqnarray}  
In particular, for all bounded measurable functions $\eta$ on $\R^d$,
\begin{equation} \int_{\R^d}\eta(\nabla\phi(v))F_0{\rm d}v = 
\int_{\R^d}\eta(w)F_1{\rm d}w\ ,\label{dfg} \end{equation}
and $\nabla \phi$ is the {\it  unique} gradient of 
 a convex function defined on 
the convex hull of the support of $F_0$ so that \pref{dfg} holds for all such 
$\eta$.

Recall that for any convex function $\psi$ on $\R^d$,  $\psi^*$ denotes its 
Legendre transform; i.e., the dual convex function, which is defined through
\begin{equation} \psi^*(w) = \sup_{v\in \R^d}\{\ w\cdot v - \psi(v)\ \}\ .\label{legendre} \end{equation}
The convex functions $\psi$ arising as optimizers in \pref{apform} have the 
further property that $(\psi^*)^* = \psi$.
Being convex, both $\psi$ and $\psi^*$ are locally Lipschitz and differentiable 
on the complement of a set of 
Hausdorff dimension $d-1$. (It is for this reason that we work with densities 
instead of measures; $\nabla
\psi\# \mu$ might not be well defined if $\mu$ charged sets Hausdorff dimension 
$d-1$.)
In our quotation of Brenier's result concerning in \pref{apform}, 
the statement that the convex functions $\psi$ and
$\phi$ in \pref{apform} are a dual pair simply means that
$\phi = \psi^*$ and $\psi = \phi^*$.
It follows from \pref{legendre} that
$\nabla\psi$ and $\nabla\psi^*$ are inverse transformations in that
\begin{equation} \nabla\psi(\nabla\psi^*(w)) = w\qquad{\rm and}
\qquad\nabla\psi^*(\nabla\psi(v)) = v\ \label{inversee} \end{equation}
for $F_1(w){\rm d}w$ almost every $w$ and $F_0(v){\rm d}v$ almost every $v$ 
respectively.

Given a map $T:\R^d\rightarrow \R^d$ and $F\in{\cal P}$, define
$T\# F\in {\cal P}$ by
$$\int_{\R^d}\eta(v)\left(T\# F(v)\right){\rm d}v =
\int_{\R^d}\eta(T(v))F(v){\rm d}v$$ 
for all test functions $\eta$ on $\R^d$. Then we can express \pref{dfg} more 
briefly by writing
$\nabla \phi \# F_0 = F_1$. The uniqueness of the gradient of the convex 
potential $\phi$ is very useful for computing $W_2^2(F_0,F_1)$ since if one can 
find {\it some} convex function $\tilde\phi$ such that
$\nabla \tilde\phi\# F_0 = F_1$, then $ \tilde\phi$ is the potential for 
the minimizing map
and
\begin{equation} W_2^2(F_0,F_1) = 
\int_{\R^d}{1\over 2}|v- \nabla\tilde\phi(v)|^2F_0(v){\rm d}v\ .\label{flick} \end{equation}

Now it is easy to determine the geodesics. These are given in terms of a natural 
interpolation between two densities $F_0$ and $F_1$ that was introduced
and applied by  McCann in his thesis \cite{mcthesis} and in  
\cite{McCann95}.

Fix two densities $F_0$ and $F_1$ in ${\cal P}$. Let $\psi$ be the 
convex
function on $\R^d$ such that $\left(\nabla\psi\right)\# F_0 = F_1$.
Then for any $t$ with $0<t<1$, define the convex function $\psi_t$ by
\begin{equation} \psi_t (v) = (1-t){|v|^2\over 2} + t\psi(v)\label{geoA1} \end{equation}
and define the density $F_t$ by
\begin{equation} F_t = \nabla \psi_t\# F_0\ .\label{geogeo} \end{equation} 
At $t=0$, $\nabla \psi_t$ is the identity, while at $t=1$, it is $\nabla\psi$. 

Clearly for each $0\le t \le 1$, $\psi_t$ is convex, and so 
the map $\nabla \psi_t$ gives the optimal transport from $F_0$ to $F_t$. What 
map gives the optimal
transport from $F_t$ onto $F_1$? 

By definition $\nabla \psi_t\# F_0 = F_t$. It follows from \pref{inversee} that
$\nabla (\psi_t)^*\# F_t = F_0$, and therefore that
$\nabla \psi\circ \nabla (\psi_t)^* \# F_t = F_1$. It turns out that  $\nabla 
\psi\circ \nabla (\psi_t)^*$
is the {\it optimal} 
transport from $F_t$ onto $F_1$. This composition property of the optimal 
transport maps
along a McCann interpolation path provides the key to  several of the theorems 
in the next section,
and is the basis of short proofs of other known results. It is the essential 
observation made in
this section.

To see that $\nabla 
\psi\circ \nabla (\psi_t)^*$ is the optimal transport map from $F_t$ onto 
$F_1$,
it suffices to show that it is a convex function. From \pref{geoA1}, 
$\nabla \psi_t(v) = (1-t)v\break + t\nabla\psi(v)$, which is the same as
$t\nabla \psi(v) = \left(\nabla \psi_t(v) - (1-t)v\right)$.
Then by \pref{inversee}, 
\begin{equation} \nabla \psi\circ \nabla(\psi_t)^*(w) = {1\over t}\left( w - (1-t)\nabla 
(\psi_t)^*(w)\right)\ .\label{wheat3} \end{equation}
Thus, $\nabla \psi\circ \nabla(\psi_t)^*(w)$ is a gradient. There are
at least two ways to proceed from here. Assuming sufficient regularity of $\psi$ 
and $\psi^*$,
one can differentiate \pref{inversee} and see that ${\rm Hess}\, \psi(\nabla\psi^*(w)){\rm Hess}\, \psi^*(w) = I$.
That is, the Hessians of $\psi$ and $\psi^*$ are inverse to one another. Since 
${\rm Hess}\, \psi_t(v) \ge (1-t)I$, this provides an upper bound on the Hessian of 
$(\psi_t)^*$ 
which can be used to show that the right side of \pref{wheat3} is the gradient of 
a convex function.
This can be made rigorous in our setting, but the argument is somewhat 
technical, and involves the definition of the 
Hessian in the sense of Alexandroff.

There is a much simpler way to proceed. As McCann showed \cite{mcthesis}, if 
$\tilde F_t$ 
is the path one gets interpolating between $F_0$ and $F_1$ but starting at 
$F_1$, then 
$F_t = \tilde F_{1-t}$. So $\nabla\left((\psi^*)_{1-t}\right)^*$ is the optimal 
transport map
from $F_t$ onto $F_1$. This tells us which convex function should have   
$\nabla \psi\circ \nabla(\psi_t)^*(w)$ as its gradient, and this is easily 
checked using the mini-max theorem.

\proclaimtitle{Interpolation and Legendre transforms} 
\proclaim{Lemma} 
 Let $\psi$ be a convex function such that $\psi = \psi^{**}$. Then
by the interpolation in {\rm \pref{geoA1},}
\begin{equation} \left((\psi^*)_{1-t}\right)^*(w) = {1\over t}\left({|w|^2\over 2} - 
(1-t)(\psi_t)^*(w)\right)\ .\label{colt1} \end{equation}
\endproclaim

\demo{Proof} Calculating, with use of the the mini-max theorem, one 
has
\begin{eqnarray*}
\left((\psi^*)_{1-t}\right)^*(w) &= &
\sup_{z}\left\{ z\cdot w - \left(t{|z|^2\over 2} + 
(1-t)\psi^*(z)\right)\right\}\\[4pt]
&=&\sup_{z}\left\{ z\cdot w - t{|z|^2\over 2} - (1-t)\sup_v\left\{v\cdot z - 
\psi(v)\right\}\right\}\\[4pt]
&=&\sup_z\inf_v\left\{z\cdot(w-(1-t)v)- t{|z|^2\over 2} + (1-t)\psi(v)\right\}\\[4pt]
&=&\inf_v\sup_z\left\{z\cdot(w-(1-t)v)- t{|z|^2\over 2} + (1-t)\psi(v)\right\}\\[4pt]
&=&{1\over t}\left({|w|^2\over 2} - (1-t)(\psi_t)^*(w)\right)\ .\\
\noalign{\vskip-36pt}
\end{eqnarray*}
\enddemo
\vglue12pt

As an immediate consequence, 
\begin{equation} \nabla  \left((\psi^*)_{1-t}\right)^* = \nabla \psi\circ \nabla 
(\psi_t)^*\label{corn45} \end{equation}
is the optimal transport from $F_t$ to $F_1$. This also implies that 
$\nabla\psi_t\#F_0 = \nabla(\psi^*)_{1-t}\#F_1$, as shown by McCann
in \cite{mcthesis} using a ``cyclic monotonicity'' argument. Lemma 2.1 leads 
to 
a simple proof of another result of McCann, 
again from  \cite{mcthesis}:

\proclaimtitle{Geodesics for the $2$-Wasserstein metric}
\proclaim{Theorem} 
Fix two densities $F_0$ and $F_1$ in ${\cal P}$. Let $\psi$ be the 
convex
function on $\R^d$ such that $\left(\nabla\psi\right)\# F_0 = F_1$.
Then for any $t$ with $0<t<1${\rm ,} define the convex function $\psi_t$ by
{\rm \pref{geoA1}}
and define the density $F_t$ by {\rm \pref{geogeo}.}
Then for all $0<t<1${\rm ,}
\begin{equation} W_2(F_0,F_t) = tW_2(F_0,F_1)\quad{\it and}\quad W_2(F_t,F_1) = 
(1-t)W_2(F_0,F_1)\quad \label{flock} \end{equation}
and $t\mapsto F_t$ is the unique path from $F_0$ to $F_1$ for the $2$\/{\rm -}\/Wasserstein 
metric that has
this property. In particular{\rm ,} there is exactly one geodesic for the\break 
$2$\/{\rm -}\/Wasserstein metric
connecting any two densities in
${\cal P}$.
\endproclaim
 
\demo{Proof} It
follows from \pref{flick} that
\begin{eqnarray*}
W_2^2(F_0,F_t) &= &
{1\over 2}\int_{\R^d}\left|v - \left((1-t)v + t\nabla\psi(v)\right)\right|^2
F_0(v){\rm d}v\\[4pt] &=&
t^2{1\over 2}\int_{\R^d}|v - \nabla\psi(v)|^2
F_0(v){\rm d}v = t^2W_2^2(F_0,F_1)\ .\end{eqnarray*}

Next, since $\nabla  \left((\psi^*)_{1-t}\right)^*$ is the optimal transport 
from $F_t$ to $F_1$,
by \pref{colt1},
\begin{eqnarray*}
W_2^2(F_t,F_1) &= &
{1\over 2}
\int_{\R^d}\left| w - {1\over t}\left(w -(1-t)\nabla 
(\psi_t)^*(w)\right)\right|^2
F_t(v){\rm d}v\\[4pt]
&= &\left({1-t\over t}\right)^2{1\over 2}\int_{\R^d}|v - \nabla\psi_t(v)|^2
F_0(v){\rm d}v = (1-t)^2W_2^2(F_0,F_1)\ .\end{eqnarray*}
Together, the last two computations give us \pref{flock}.

The uniqueness follows from a strict convexity property of the distance: 
For any probability density $G_0$, 
the function $G\mapsto  W_2^2(G_0,G)$ is strictly convex on ${\cal P}$
in that for any pair $G_1$, $G_2$ in ${\cal P}$ and any $t$ with
$0 < t < 1$,
\begin{equation} W_2^2(G_0,(1-t)G_1 + tG_2) \le (1-t)W_2^2(G_0,G_1) +
tW_2^2(G_0,G_2)\hskip.4in \label{apstrictcon} \end{equation}
and there is equality if and only if $G_1 = G_2$. 
 This follows easily from the uniqueness of the optimal coupling specified in 
 \pref{apform}; nontrivial
convex combinations of such couplings are not of the form \pref{apform}, and 
therefore cannot be optimal. 

Now suppose that there are two geodesics $t\mapsto F_t$ and $t\mapsto \tilde 
F_t$. Pick some
$t_0$ with $F_{t_0} \ne \tilde F_{t_0}$. Then the path consisting of a geodesic 
from
$F_0$ to $( F_{t_0} +\tilde F_{t_0})/2$, and from there onto $F_1$ would have a 
strictly shorter length
than the geodesic from $F_0$ to $F_1$, which cannot be. \enddemo

To obtain an Eulerian description of these geodesics, let $f$ be any smooth 
function on $\R^d$, and compute:
\begin{eqnarray}
&&\label{rice1} \\
{{\rm d} \over {\rm d}t}\int_{R^d}f(v)F_t(v){\rm d}v &=&
{{\rm d} \over {\rm d}t}\int_{R^d}f(\nabla\psi_t(v))F_0(v){\rm d}v\nonumber\\
&=& \int_{R^d} \nabla f(\nabla\psi_t(v))\left[v- \nabla \psi(v)\right]F_0(v){\rm 
d}v\nonumber \\
&= &\int_{R^d} \nabla f(w)\left[\nabla (\psi_t)^*(w)- \nabla \psi(\nabla 
(\psi_t)^*(w))\right]F_t(w){\rm d}w\nonumber \\
&= &\int_{R^d} \nabla f(w)\left[{w - \nabla (\psi_t)^*(w)\over 
t}\right]F_t(w){\rm 
d}w\ .\nonumber
\end{eqnarray}
In other words, when $F_t$ is defined in terms of $F_0$ and $\psi$ as in 
\pref{geoA1} and \pref{geogeo},  
$F_t$ is a weak solution to
\begin{equation} {\partial\over \partial t}F_t(w) + \nabla\cdot\left(W(w,t)F_t(w)\right)= 
0\label{wheat4} \end{equation}
where, according to Lemma 2.1,
\begin{equation} W(w,t) = {w - \nabla (\psi_t)^*(w)\over t} = \nabla\left({|w|^2\over 2t} - 
{1\over t}(\psi_t)^*(w)\right)\ .\label{wheat5} \end{equation}
In light of the first two equalities in \pref{rice1}, 
\begin{equation} W(w,0) = \nabla\left({|w|^2\over 2} - \psi(w)\right) = w -\nabla \psi(w)\ 
.\label{wheat6} \end{equation}
This gradient vector field can be viewed as giving the ``tangent direction'' 
to the geodesic $t\mapsto F_t$ at $t=0$.

We would like to identify some subspace of the space of gradient vector fields
as the tangent space $T_{F_0}$  
to ${\cal P}$ at $F_0$. Towards this end we ask: Given a smooth, rapidly 
decaying function $\eta$ on $\R^d$, 
is there a geodesic $t\mapsto F_t$
passing through $F_0$ at $t=0$ so that, in the weak sense,
\begin{equation} \left({\partial \over \partial t}F_t + \nabla \cdot \left(\nabla \eta 
F_t\right)\right)\bigg|_{t=0} = 0\ .\label{rice2} \end{equation}
The next theorem says that this is the case, and provides us with a geodesic 
that \pref{rice2} holds with 
$\eta$ sufficiently small. But
then by changing the time parametrization, we obtain a geodesic, possibly quite 
short, that has any multiple of $\nabla\eta$ as 
its initial ``tangent vector''.

 \proclaimtitle{Tangents to geodesics} 
\proclaim{Theorem}Let $\eta$ be any 
smooth{\rm ,} 
rapidly decaying function $\eta$ on $\R^d$ such that for all $v${\rm ,}
\begin{equation} \psi(v) = {|v|^2\over 2} +\eta(v)\label{wheat7} \end{equation}
is strictly convex. For any density $F_0$ in ${\cal P}${\rm ,}
and $t$ with $0\le t\le 1${\rm ,}  define 
\begin{equation} \nabla \psi_t(v) = (1-t)v + t\nabla  \psi(v) = v + t\nabla\eta(v)\ .
\label{gradform} \end{equation}
Then for all $t$ with $0\le t \le 1${\rm ,} $F_t = \nabla \psi_t\# F_0$ is absolutely 
continuous{\rm ,} and
is a weak solution of
\begin{equation} {\partial \over \partial t}F_t(v) + \nabla\cdot\left(
\nabla\eta_t(v)F_t(v)\right) =0\ ,\label{tangent}
\end{equation}
where
\begin{equation} \eta_t(v) = {1\over t}\left({|v|^2\over 2} - 
(\psi_t)^*(v)\right).\label{etadef3} \end{equation}
Moreover{\rm ,}
\begin{equation} \nabla \eta_t(v) =  \nabla\eta(v) - {t\over 2}\nabla|\nabla\eta(v)|^2 + 
t^2\nabla R_t(v)\ ,\label{etatder} \end{equation}
where the remainder term $\nabla R_t(v)$ satisfies $\|\nabla R_t\|_\infty \le 
\|{\rm Hess}\, (\eta)\|_\infty^2$
uniformly in $t$.
\endproclaim

\demo{Proof} First, the fact that $\nabla\psi_t\#F_0$ is absolutely 
continuous follows from the fact that $\nabla (\psi_t)^*$
is Lipschitz. Formulas \pref{tangent} and \pref{etadef3} follow directly from 
\pref{wheat4} and \pref{wheat5}.

To obtain \pref{etatder}, use \pref{inversee} to see that $\nabla (\psi_t)^*(v) = 
\Phi(\nabla (\psi_t)^*(v))$
where $\Phi(w) = v - t\nabla\eta(w)$. Iterating this fixed point equation three 
times yields \pref{etatder}.
\enddemo

In light of Theorems 2.2 and  2.3, we now know that every geodesic $t\mapsto 
F_t$ through $F_0$ at $t=0$ satisfies \pref{rice2}, and conversely, for every 
smooth rapidly decaying gradient vector field, there is a geodesic $t\mapsto 
F_t$ through $F_0$ at $t=0$ satisfying \pref{rice2} for that function $\eta$. 
Moreover,  along this geodesic
\begin{equation} W_2^2(F_0,F_t) = \int_0^t \left(\int_{\R^d}|\nabla\eta_s(v)|^2F_s(v){\rm 
d}v\right){\rm d}s
= t\int_{\R^d}|\nabla\eta(v)|^2F_0(v){\rm d}v\ ,\label{rice6} \end{equation}
where $\eta_s$ is related to $\eta$ as in Theorem 2.3.

Furthermore if $t\mapsto F_t$ is a path in ${\cal P}$ satisfying  \pref{rice2}
for some gradient vector field $\nabla\eta$, then this vector field is unique. 
For suppose that
$t\mapsto F_t$ also satisfies
\begin{equation} \left({\partial \over \partial t}F_t + 
\nabla \cdot \left(\nabla \xi F_t\right)\right)\bigg|_{t=0} = 0\ .\label{rice12} \end{equation}
Then, 
$\nabla \cdot \left(\nabla (\eta - \xi) F_0 \right)= 0$. Integrating against 
$\eta- \xi$, we obtain that
$$\int_{\R^d}|\nabla\eta - \nabla\xi|^2F_0(v){\rm d}v= 0\ .$$
Careful consideration of this well-known argument, inserting a cut-off function 
before integrating by parts,
reveals that all it requires is that both $\nabla\eta$ and $\nabla\xi$ are 
square integrable with respect to $F_0$. This justifies the identification of 
the tangent vector $\partial F/\partial t$ 
with $\nabla\eta$ when \pref{rice2} holds and
$\nabla\eta$ is square integrable with respect to $F_0$.

This identifies the ``tangent vector'' $\partial F_t/ \partial t$ with $\nabla 
\eta$, and gives us the 
Riemannian metric, first introduced by Otto \cite{Ottop},
\begin{equation} g\left({\partial F\over \partial t},{\partial F \over \partial t}\right) 
={1\over 2}\int|\nabla\eta(v)|^2F_0(v){\rm d}v \ .\label{metric} \end{equation}
By \pref{rice6}, the distance on ${\cal P}$ induced by this metric is the 
$2$-Wasserstein distance.

Interestingly, Theorem 2.2 provides a global description of the geodesics 
without having to
first determine and study the Riemannian metric. Theorem~2.3 gives an Eulerian 
characterization of the geodesics
which provides a complement to McCann's original Lagrangian characterization. 
Another Eulerian
analysis of the geodesics in terms of the Hamilton-Jacobi equation seems to be 
folklore in the subject.
A clear account can be found in recent lecture notes of Villani 
\cite{Villani}.

We now turn to the notion of convexity on ${\cal P}$ with respect to the\break 
$2$-Wasserstein metric. A functional
$\Phi$ on ${\cal P}$ is said to be {\it displacement convex at} $F_0$ in case 
$t\mapsto \Phi(F_t)$
is convex on some neighborhood of $0$ for all geodesics $t\mapsto F_t$ passing 
through $F_0$ at $t=0$.
A functional $\Phi$ on ${\cal P}$ is said to be {\it displacement convex} if it 
is displacement convex at all points $F_0$ of 
${\cal P}$.

If moreover $t\mapsto \Phi(F_t)$ is twice differentiable, we can check for 
displacement convexity by computing
the Hessian:
\begin{equation} {\rm Hess}\, \Phi(F_0)\langle \nabla\eta,\nabla\eta\rangle =
{{\rm d}^2\over {\rm d}t^2}\Phi(F_t)\biggl|_{t=0}\ ,\label{hess} \end{equation}
where $\nabla\eta$ is the tangent to the geodesic at $t=0$.

\proclaimtitle{Displacement convexity}
\proclaim{Theorem}  If 
the functional $\Phi$ on  ${\cal P}$
is given by
\begin{equation} \Phi(F) = \int_{\R^d}g(F(v)){\rm d}v  \label{gform} \end{equation}
where $g$ is a twice differentiable convex  function on $\R_+${\rm ,}
then $\Phi$ is displacement convex if  
\begin{equation} tg'(t) - g(t) \ge 0\qquad{\it and}\qquad t^2g''(t) -tg'(t) + g(t)\ge 0   
\label{gcond} \end{equation}
for all $t> 0${\rm ,} where the primes denote derivatives. 
\endproclaim

\demo{Proof} We check for convexity at a density $F_0$ in the domain of 
$\Phi$. 
By a standard mollification, we can find a sequence of smooth densities 
${F_0^{(n)}}$ 
with $\lim_{n\to\infty}F_0^{(n)} = F_0$ and  $\lim_{n\to\infty}\Phi(F_0^{(n)}) = 
\Phi(F_0)$.
Fix any smooth rapidly decaying function $\eta$, such that (taking a small 
multiple if need be)
$|v|^2 + \eta(v)$ is strictly convex. Then with $\nabla \psi_t$ defined as in 
\pref{gradform},
$$t\mapsto \nabla \psi_t \#F_0^{(n)}  = F_t^{(n)}$$
gives a geodesic passing through $F_0^{(n)}$ at $t=0$ with the tangent direction 
$\nabla\eta$, and
defined for $0\le t \le 1$ uniformly in $n$. Also,
$\lim_{n\to\infty}\Phi(F_t^{(n)}) = \Phi(F_t)$ for all such $t$. Therefore, it 
suffices to show that
for each $n$, $t\mapsto  \Phi(F_t^{(n)})$ is convex. In other words, we may 
assume that $F_0$ is smooth.
Then so is each $F_t$, since $F_t(w) = F_0(\nabla(\psi_t)^*(w)){\rm 
det}\left({\rm Hess}\, (\psi_t)^*)(w)\right)$
is a composition of smooth functions.
We may now check convexity by differentiating.

By \pref{tangent},
\begin{eqnarray*}
{{\rm d}\over {\rm d}t}\int_{\R^d}g(F_t(v)){\rm d}v &= &
-\int_{\R^d}g'(F_t(v))\nabla\cdot\left(\nabla \eta_t(v)F_t(v)\right){\rm d}v \\
&=&\int_{\R^d}\left(g''(F_t(v))\nabla F_t(v)\right)\cdot\left(\nabla 
\eta_t(v)F_t(v)\right){\rm d}v.\end{eqnarray*}
Defining 
$h(t) = tg'(t) - g(t)$ so that $h'(t) = tg''(t)$,
one has from \pref{tangent} that
\begin{equation}{{\rm d}\over {\rm d}t}\Phi(F_t)= 
\int_{\R^d}\nabla h(F_t(v))\cdot\nabla\eta_t(v){\rm d}v\ .\label{firstder} \end{equation}
To differentiate a second time, use \pref{etatder} to obtain
$${{\rm d}^2\over {\rm d}t^2}\Phi(F_t)\biggl|_{t=0} = \int_{\R^d}\nabla h(F_0)
\cdot \nabla\left(-{1\over 2}|\nabla \eta|^2\right){\rm d}v -
\int_{\R^d} {\partial\over \partial t}h(F_t)\biggl|_{t=0}
\left(\Delta \eta\right){\rm d}v\ .$$
But
$$ {\partial\over \partial t}h(F_t)\biggl|_{t=0}=  - F_0^2g''(F_0)\left(\Delta 
\eta\right) - 
\nabla h(F_0)\cdot\nabla\eta$$
and hence
\begin{eqnarray}
&&\label{poshes}\\
&&\hskip-12pt {{\rm d}^2\over {\rm d}t^2}\Phi(F_t)\biggl|_{t=0} \nonumber\\[4pt]
&&\qquad =\int_{\R^d}\nabla h(F_0)
\cdot \left(-{1\over 2}\nabla|\nabla \eta|^2 + 
\left(\Delta\eta\right)\nabla\eta\right){\rm d}v +
\int_{\R^d} F_0^2g''(F_0)\left(\Delta \eta\right)^2{\rm d}v \nonumber\\
&&\qquad =\int_{\R^d} h(F_0)\|{\rm Hess}\, \eta\|^2{\rm d}v +
\int_{\R^d} \left(F_0^2g''(F_0) - h(F_0)\right)
\left(\Delta \eta\right)^2{\rm d}v\ .\nonumber   \end{eqnarray}
Here, $\|{\rm Hess}\, \eta\|^2$ denotes the square of the Hilbert-Schmidt norm of
the Hessian of $\eta$. This quantity is positive whenever
$h(F) = Fg'(F) - g(F)$ and $F^2g''(F) - h(F) = F^2g''(F) - Fg'(F) + g(F)$ are 
positive.
\enddemo

The case of greatest interest here is the entropy functional  $S(F)$, defined in 
\pref{xent}.
In this case,
$g(t) = t\ln t$, so that  $tg'(t) -g(t) = t$ and $tg''(t) - tg'(t) +g(t) = 0$.
Hence from \pref{poshes},
\begin{equation} {{\rm d}^2\over {\rm d}t^2}S( F_t)\biggl|_{t=0} =
\int_{\R^d} \|{\rm Hess}\, \eta\|^2F_0(v){\rm d}v\ .\label{ent1} \end{equation}

This shows that the entropy is convex, 
as proved  in \cite{Ottop}, though not strictly convex. Consider the 
following example\footnote{We thank the referee for this example, which has clarified the formulation of 
Corollary 2.5 below.}
in one dimension:
Let 
$$\psi(v) = {|v|^2\over 2} + |v|.$$
For any $F_0$, 
define $F_t = \nabla \psi_t$ and then it is easy to see that
\begin{equation} F_t(v) = 1_{\{v <-t\}}F_0(v+t) + 1_{\{v > t\}}F_0(v-t)\ .\label{colt57} \end{equation}
The geodesic $t\mapsto F_t$ can be continued indefinitely for positive $t$, but 
unless $F_0$ vanishes in some strip $-\varepsilon < v <\varepsilon$, 
it cannot be continued at all for negative $t$.
With $F_t$ defined as in \pref{colt57}, $S(F_t) = S(F_0)$ for all $t$.

There are however 
interesting cases in which the entropy is strictly convex along a geodesic, and 
even uniformly so:
Suppose that the ``center of mass'' $\int_{\R^d}vF_t(v){\rm d}v$ is
constant along the geodesic $t \mapsto F_t$, which means that
\begin{equation} \int_{\R^d} \nabla\eta(v)F_0(v){\rm d}v =0\label{null} \end{equation}
where as above, $\nabla\eta$ is the tangent vector generating the geodesic.

The Poincar\'e constant $\alpha(F)$ of a density $F$ in ${\cal P}$ is defined by
\begin{equation} \alpha(F) = \inf_{\varphi\in {\cal C}_0^\infty}
{\int |\nabla \varphi(v)|^2F(v){\rm d}v \over 
\int\left|\varphi(v) - \int \varphi(v)F(v){\rm d}v\right|^2F(v){\rm d}v}\ 
.\label{poincare} \end{equation}
Thus, when \pref{null} holds, with $\varphi = {\partial \eta/\partial v_i}$ for 
$i = 1\dots d$
 we take the sum,
yielding
\begin{equation} \int_{\R^d} \|{\rm Hess}\, \eta\|^2F_0(v){\rm d}v \ge 
\alpha(F_0)\int_{\R^d} |\nabla \eta(v)|^2F_0(v){\rm d}v\ ,\label{colt56} \end{equation}
which provides a lower bound to the right side of \pref{ent1} in terms of the 
Riemannian metric.

Now consider a  ``smooth'' geodesic through a smooth density $F_0$, as in the 
previous proof, 
and such that \pref{null} is satisfied. Then by (2.31) and \pref{colt56},
for any $t$ and $h>0$ such that $F_{t-h}$ and $F_{t+h}$ are both on the 
geodesic,
$${1\over h^2}\left(S(F_{t+h})+S(F_{t-h}) -2S(F_t)\right) \ge \alpha(F_t)\
\int_{\R^d} |\nabla \eta(v)|^2F_0(v){\rm d}v\ .$$
If the geodesic is parametrized by arclength, 
then the last factor on the right is one.

Summarizing the last paragraphs, we have the following corollary:

\proclaimtitle{Strict convexity of entropy} 
\proclaim{{C}orollary} Consider a 
geodesic 
$s\mapsto F_s$ parametrized by arc length $s${\rm ,} and defined for some interval 
$a< s< b$
such that
$s\mapsto \int vF_s(v){\rm d}v$ is constant{\rm ,} and such that each $F_s$ is bounded 
and continuously differentiable. Then for all $s$ and  $h$ so that $a < 
s-h,s+h<b${\rm ,}
\begin{equation} S(F_{s+h})+S(F_{s-h}) -2S(F_s)
\ge h^2\alpha(F_s)\ ,
\label{strict} \end{equation}
where $\alpha(F_s)$ is the Poincar{\rm \'{\it e}} constant of the density $F_s$.
\endproclaim

(Notice that for the geodesic \pref{colt57}, $\alpha(F_t) = 0$ for all $t>0$, as 
long as $F_0$ has positive mass on both sides of the origin, in addition to the 
fact that
$F_t$ will not in general be smooth.)

We remark that Caffarelli has recently shown \cite{Caff00} that if $F_0$ is a 
Gaussian density,
and $F_1 = e^{-V}F_0$ where $V$ is convex, then there is an {\it upper} bound on 
the Hessian of 
the potential $\psi$ for which $\nabla\psi\# F_0  =
F_1$. This upper bound is inherited by $\psi_t$ for all $t$. Since as Caffarelli 
shows, an upper bound on
the Hessian of $\psi$ and a lower bound on the Poincar\'e constant for $F_0$ 
imply a lower bound on
the Poincar\'e constant of $F_t$, one obtains a {\it uniform} 
lower bound on the 
Poincar\'e constant for $F_t$,
$0<t<1$. Hence $S(F_t)$ is uniformly strictly convex along such a geodesic.

\section{Geometry of the constraint manifold}

Let $u \in \R^d$ and $\theta > 0$ be given.
Consider the subset
${\cal E}_{u,\theta}$ of ${\cal P}$
specified by
\begin{equation} {\cal E}_{u,\theta} = \left\{F\in{\cal P}\ \biggl|\ {1\over 
d}\int_{\R^d}|v-u|^2 
F(v){\rm d}v = \theta
\quad{\rm and}\quad \int_{\R^d}vF(v){\rm d}v = u\ \right\}.\enspace\label{mdef} \end{equation}
This is the set of all probability
densities with a mean $u$ and a variance $d\theta$. 
We will often
write ${\cal E}$ in place of ${\cal E}_{u,\theta}$ when 
$u$ and $\theta$ are clear from the context or simply irrelevant.

We give a fairly complete description of the geometry of ${\cal E}$, both 
locally and globally.
In particular, we  obtain a closed form expression for the distance between any 
two points on ${\cal E}$
in the metric induced by the $2$-Wasserstein metric, and  a global description of 
the geodesics in
${\cal E}$.

Notice that 
\begin{equation} {\cal E}_{u,\theta} \subset \left\{ F\ | \ W_2^2(F,\delta_u) = {d\theta\over 
2}\ \right\}\label{UU45} \end{equation}
where $\delta_u$ is the unit mass at $u$. This is quite clear from the 
transport point of view: If our target distribution is a point mass, 
there are no choices to make; everything is simply transported to the point
$u$. Hence   ${\cal E}_{u,\theta}$ is 
a part of a  sphere in the 2-Wasserstein metric, 
 centered on $\delta_u$, and with a radius of $\sqrt{d\theta/2}$.

Our first theorem shows that for any $F_0$ in ${\cal P}$, there is a unique 
closest $F$ in ${\cal E}$, and this is obtained by dilatation and translation. 
This is the first of two 
related variational
problems solved in this section.

 \proclaimtitle{Projection onto ${\cal E}$}
\proclaim{Theorem}  Let $F_0$ be any 
 probability density on 
$\R^d$ such that  
$$\int_{\R^d}v F_0(v){\rm d}v = u_0 \qquad{\it and}\qquad\int_{\R^d}|v-u_0|^2 
F_0(v){\rm d}v = d\theta_0\ .$$
Let $\theta>0$ and $u$ be given{\rm ,} and set $a = \sqrt{\theta_0/\theta}$. Then
$$\inf\left\{ W_2^2(G,F_0)\ |\ G \in {\cal E}_{\theta,u}\right\}$$
is attained at
$$ \tilde F(v) = a^d F_0\left(a(v-u)+u_0 \right)\ ,$$
and the minimum value is
\begin{equation} W_2^2(F_0,\tilde F) = {\left(\sqrt{\theta} - \sqrt{\theta_0}\right)^2\over 2} 
+ {|u -u_0|^2\over 2}\ .\label{fanmin} \end{equation}
\endproclaim

\demo{Proof} There is no loss of generality in fixing $u=0$ in the 
proof
since if $u_0$ is arbitrary, a translation of both $\tilde F$ and $F_0$ yields 
the general result.

Let $\phi$ be defined by $\phi(v) =
|v-u_0|^2/(2a)$ so that $\left(\nabla\phi\right)\# F_0 =  \tilde F$. 
Let $\psi(w) = a|w|^2/2+w\cdot u_0$ be the dual convex function so that
$$\phi(v) + \psi(w) \ge v\cdot w\ ,$$
and hence
\begin{equation} {1\over 2}|v-w|^2 \ge {a|v|^2 - |v-u_0|^2\over 2a} +
{(1-a)|w|^2-w\cdot u_0\over 2}\label{young} \end{equation}
for all $v$ and $w$. 

Next, given any $G$ in ${\cal E}$, let $\gamma$ be
the optimal coupling of $F_0$ and $G$ so that
$$W_2^2(F_0, G) = \int_{\R^d\times\R^d}{1\over 2}|v-w|^2\gamma({\rm 
d}v,{\rm d}w) \ .$$
Then by \pref{young},
\begin{eqnarray*}
W_2^2(F_0,G)  &\ge& 
 \left({a-1\over 2a}\right)\int_{\R^d}|v|^2
F_0(v){\rm d}v + {|u_0|^2\over 2} + 
\left({1-a\over 2}\right)\int_{\R^d}|w|^2G(w){\rm d}w \\
&=& {(a-1)^2d\theta\over 2}+ {|u_0|^2\over 2}\ .\end{eqnarray*}
On the other hand, since $\left(\nabla\phi\right)\# F_0 =  \tilde F$,
\begin{eqnarray*}
W_2^2(F_0,\tilde F) &= &
\int_{\R^d}{1\over 2}|v - \nabla\phi(v)|^2F_0(v){\rm d}v \\
&=& \left({1\over a}-1\right)^2\int_{\R^d}{1\over 2}|v-u_0|^2
F_0(v){\rm d}v + {|u_0|^2\over 2} \\ 
&=& {(a-1)^2d\theta\over 2} + {|u_0|^2\over 2}\ .\\
\noalign{\vskip-36pt}\end{eqnarray*}
\enddemo
\vglue12pt

\demo{{R}emark {\rm (Exact solution for the {\rm JKO} time discretization of the
heat equation for Gaussian initial data)}}
Theorem 3.1 allows us to solve exactly  the Jordan-Kinderlehrer-Otto time 
discretization of the
heat equation for Gaussian initial data. Take as initial data $F_0(v) = (4\pi 
t_0)^{-d/2}e^{-|v|^2/4t_0}$. We can now find
$\inf\{ W_2^2(F,F_0) + hS(F) \}$ in two steps. First, consider
\begin{equation} 
\inf\{ W_2^2(F,F_0) + hS(F)\ |\ F \in {\cal E}_{0, 2t d} \}.\label{heatinft}
\end{equation}
Now on ${\cal E}_{0, 2t d}$, $S$ has a global minimum at $G_t = (4\pi 
t)^{-d/2}e^{-|v|^2/4t}$, as is well known. By 
Theorem 3.1,
$ W_2^2(F,F_0)$ also has a global minimum on ${\cal E}_{0, 2t d}$ at $G_t$, 
since $G_t$ is just a rescaling of $F_0$.
Therefore, by \pref{fanmin}, the infimum in \pref{heatinft} is 
$$W_2^2(G_t,F_0) + hS(G_t) = d\left(\sqrt{t} - \sqrt{t_0}\right)^2 + -h{d\over 
2}\left(\ln(4\pi  t) + 1\right)\ .$$
In the second step, we simply compute the minimizing value of $t$, which amounts 
to finding the value of $t$ that minimizes
$$\left(\sqrt{t} - \sqrt{t_0}\right)^2 - {h\over 2 }\ln t\ . $$
Simple computations lead to the value
$t = f(t_0)$ where
\begin{equation} f(s) = {1\over 2}\left(s+h+ s\sqrt{1+ {2h\over s}}\right)\ .\label{iterf} \end{equation}
Note that $t_0 < f(t_0) < t_0+h$, but $f(t_0) = t_0+h+O(h^2)$. If we then 
inductively define $t_n = f(t_{n-1})$, we see that the
exact solution of the Jordan-Kinderlehrer-Otto time discretization of the
heat equation is given at time step $n$ by $F_n = (4\pi 
t_n)^{-d/2}e^{-|v|^2/4t_n}$
where $t_n = t_0+ nh + {\cal O}(h^2)$. Note that in the discrete time 
approximation, 
the
variance increases more slowly than in continuous time, since the ${\cal 
O}(h^2)$ term is negative,
though of course the difference in the rates vanishes as $h$ tends to zero.

Returning to the main focus of this section, fix two densities $F_0$ and $F_1$ 
in ${\cal E}$. Let $\psi$ be the 
convex
function on $\R^d$ such that $\left(\nabla\psi\right)\# F_0 = F_1$.
Then by Theorem 2.2, the geodesic that runs from $F_0$ to $F_1$ through
the ambient space ${\cal P}$ is given by
$$F_t = \left((1-t)v + t\nabla\psi\right)\# F_0\ .$$
Thinking of ${\cal E}$ as a subset of a sphere, and this geodesic as the chord 
connecting two points on the sphere,
we refer to it as the  {\it chordal geodesic} $F_0$ to $F_1$.

 \proclaimtitle{Variance along a chordal geodesic}
\proclaim{Lemma} Let $F_0$ and 
$F_1$ be any two densities
in ${\cal E}$. Let $t\mapsto F_t$ be the chordal geodesic joining them. Then for 
all $t$ with $0\le t \le 1${\rm ,}
\begin{eqnarray}
{1\over 2}\int_{\R^d}|v-u|^2F_t(v){\rm d}v &=&
{d\theta\over 2} \left[1 - 4t(1-t){W_2^2(F_0,F_1)\over 2d\theta}\right]\label{chordalvar}\\
&=& R_\theta^2\left[1 - t(1-t){W_2^2(F_0,F_1)\over R^2_\theta}\right]\nonumber
 \end{eqnarray}
where $R_\theta = \sqrt{d\theta/2}$.
\endproclaim
 
\demo{Proof} Notice first that with $F_1 = \nabla \psi\#F_0$, we have 
from Theorem 2.2 that
\begin{eqnarray}
&&\int_{\R^d}{1\over 2}|v-u|^2F_t(v){\rm d}v\label{smaller} \\
&&\quad\quad = 
\int_{\R^d}{1\over 2}|\left((1-t)v + t\nabla\psi(v)\right) -u|^2F_0(v){\rm 
d}v\nonumber
\\
&&\quad\quad = \int_{\R^d}{1\over 2}|(1-t)\left(v-u\right) + t\left(\nabla\psi(v) 
-u\right)|^2F_0(v){\rm d}v\nonumber\\ 
&&\quad\quad = (1-t)^2\int_{\R^d}{1\over 2}|v-u|^2F_0(v){\rm d}v +
t^2\int_{\R^d}{1\over 2}|w-u|^2F_1(y){\rm d}v\nonumber\\  
&&\quad\qquad +t(1-t)\int_{\R^d}(v-u)\cdot \left(\nabla\psi(v)- u\right)F_0(v){\rm d}v\nonumber\\ 
&&\quad\quad = {d\theta\over 2} (1-t)^2 + {d\theta\over 2} t^2 + t(1-t)\int_{\R^d}(v-u)\cdot 
\left(\nabla\psi(v)-
u\right)F_0(v){\rm d}v.\nonumber \end{eqnarray}
Next,
\begin{eqnarray*}
W_2^2(F_0,F_1) &= &
{1\over 2}\int_{\R^d}|v - \nabla\psi(v)|^2F_0(v){\rm d}v \\
&=&{1\over 2}\int_{\R^d}|v-u|^2F_0(v){\rm d}v +
{1\over 2}\int_{\R^d}|\nabla\psi(v)-u|^2F_0(v){\rm d}v \\
&&-\int_{\R^d}(v-u)\cdot \left(\nabla\psi(v)- u\right)F_0(v){\rm d}v \\
&=&d\theta - \int_{\R^d}(v-u)\cdot \left(\nabla\psi(v)- u\right)F_0(v){\rm 
d}v\end{eqnarray*}
by the definition of ${\cal E}$, and hence
\begin{equation} \int_{\R^d}(v-u)\cdot \left(\nabla\psi(v)- u\right)F_0(v){\rm d}v = 
d\theta - W_2^2(F_0,F_1)\ .\label{vdoty} \end{equation}
Combining   \pref{vdoty} and \pref{smaller}, one has the result. \enddemo

We note that since  
$\int_{\R^d}(v-u)F_0(v){\rm d}v =0$, 
$$\int_{\R^d}(v-u)\cdot \left(\nabla\psi(v)- u\right)F_0(v){\rm d}v =
\int_{\R^d}(v-u)\cdot \left(\nabla\psi(v)- \nabla\psi(u)\right)F_0(v){\rm d}v 
\ge 0$$
by the convexity of $\psi$. It follows from this and \pref{vdoty} that
\begin{equation}   W_2^2(F_0,F_1) \le d\theta = 2R_\theta^2\ ,\label{lowvarbnd} \end{equation}
where $R_\theta = \sqrt{d\theta/2}$ is the radius of ${\cal E}$ as in 
\pref{UU45}.
Hence the variance in \pref{chordalvar} is never smaller than $R_\theta^2$.

The next result is the second of the variational problems solved in this 
section, and is the key to the determination of
the geodesics in ${\cal E}$.

\proclaimtitle{Midpoint theorem}
\proclaim{Theorem}  Let $F_0$ and $F_1$ be any 
two densities in ${\cal E}$.
Then
\begin{equation} \inf_{G\in {\cal E}}\left\{ W_2^2(F_0,G) + W_2^2(G,F_1)\right\}\label{midp} \end{equation}
is attained uniquely at $a^dF_{1/2}(a(v-u)+u)$ where $F_{1/2}$ is the midpoint 
of the chordal geodesic{\rm ,}
and $a$ is chosen to rescale the midpoint onto ${\cal E}${\rm ;} i.e.{\rm ,}
\begin{equation} a = \sqrt{1 - {W_2^2(F_0,F_1)\over 2d\theta}} = \sqrt{1 - {W_2^2(F_0,F_1)\over 
(2R_\theta)^2}}\ ,\label{aform} \end{equation}
where $R_\theta = \sqrt{d\theta/2}$ is the radius of ${\cal E}$ as in 
{\rm \pref{UU45}.}
Moreover{\rm ,} the minimal value attained in {\rm \pref{midp}} is 
$f\left(W_2^2(F_0,F_1)\right)$
where
\begin{equation} f(x) = 2d\theta\left(1 - \sqrt{1-x/(2d\theta)}\right)\ .\label{gonzo} \end{equation}
The function $f$ is convex and increasing on $[0,2d \theta]$.\endproclaim

Before giving the proof itself, we first consider some formal arguments that 
serve to identify the minimizer
and motivate the proof.

Let $\Phi(G)$ denote the functional being minimized in \pref{midp}. This 
functional
is strictly convex with respect to the usual convex structure on  ${\cal E}$; 
that is,
for all $\lambda$ with $0 < \lambda < 1$, and all $G_0$ and $G_1$ in  ${\cal 
E}$, 
$$ \Phi(\lambda G_0 + (1-\lambda)G_1) \le \lambda \Phi(G_0) + 
(1-\lambda)\Phi(G_1)$$
with equality only if $G_0 = G_1$. 
The strict convexity suggests that there is a minimizer $G_0$, and that if we 
can
find any critical point
$G$ of $\Phi$, then
$G$ is the minimizer $G_0$. 

To make variations in $G$, seeking a critical point, let $\eta$ be a smooth, 
rapidly decaying function on $\R^d$, and define the map
$T_t:\R^d\rightarrow
\R^d$ by $T_t(v) = v + t\nabla\eta(v)$. Let $G_t = T_t\#G_0$. We want the curve 
$t\mapsto G_t$ to 
be tangent to ${\cal E}$ at $t=0$, and so we require in particular that
\begin{equation} \int_{\R^d}v\cdot\nabla\eta(v)G_0(v)
{\rm d}v =0\ \label{crcon} \end{equation}
which guarantees that $\int |v|^2G(t) {\rm d}v = \int |v|^2G_0 {\rm d}v + {\cal 
O}(t^2)$.

Let $\phi$ be the convex function such that $\nabla \phi \#G_0 = F_0$,
and let $\tilde \phi$ be the convex function such that $\nabla \tilde \phi \#G_0 
= F_1$.
The variation in $\Phi(G_t)$ can be expressed in terms of $\phi$, $\tilde \phi$ 
and  $\eta$ as follows:
Formally, assuming enough regularity, we have
\begin{equation} \lim_{t\to 0^+}{\Phi(G_t)-\Phi(G_0)\over t} =
\int_{\R^d}\left(\nabla \phi(v) + \nabla\tilde\phi(v) 
-2v\right)\cdot\nabla\eta(v)G_0(v)
{\rm d}v\ .\quad\label{varf} \end{equation}
(A more precise statement and explanation are provided in Section 4 where we make 
actual use of such variations.
For the present heuristic purposes it suffices to  \pagebreak be formal.) 

Combining \pref{crcon} and \pref{varf}, we see that the formal condition for
$G_0$ to be a critical point is 
\begin{equation} \nabla \phi(v) + \nabla\tilde\phi(v) = Cv\label{critcrit} \end{equation}
for some constant $C$.

The formal argument tells us what to look for, namely a $G_0$ such that 
\pref{critcrit} holds. It is easy to see,
if $G_0$ is the midpoint  of the chordal geodesic from $F_0$ to $F_1$ 
projected onto ${\cal E}$ by rescaling as in Theorem 3.1, that $G_0$ satisfies
\pref{critcrit}. The actual proof of the theorem consists of two steps: First we 
verify the assertion just made about
$G_0$ so defined. Then we prove, using \pref{critcrit}, that $G_0$ is indeed the 
minimizer using a
duality argument very much like the one used to prove Theorem 3.1.

\demo{Proof of Theorem {\rm 3.3}} First, we may assume that $u=0$. Next, 
let $\psi$ be the convex function such that 
$\nabla \psi \#F_0 = F_1$. We may suppose initially 
that both $F_0$ and $F_1$ are strictly positive so that $\psi$ will be convex on 
all of $\R^d$.
Recall that 
$\nabla\left(\psi_{1/2}\right)^*\# F_{1/2} = F_0$, and that by \pref{corn45},
$\nabla\left((\psi^*)_{1/2}\right)^*\# F_{1/2} = F_1$. Then immediately from  
\pref{colt1} we have 
\begin{equation} \left(\psi_{1/2}\right)^*(v) + \left((\psi^*)_{1/2}\right)^*(v) = |v|^2\ 
.\label{corn73} \end{equation}

Now let $a$ be given by \pref{aform}, and define
$$\phi(v) = {1\over a}\left(\psi_{1/2}\right)^*(av)\qquad{\rm and}\qquad \tilde 
\phi 
= {1\over a}\left((\psi^*)_{1/2}\right)^*(av)\ .$$
Then, $\nabla\phi\# G_0 = F_0$ and $\nabla\tilde\phi\# G_0 = F_1$, and from 
\pref{corn73}, 
\begin{equation} \phi(v) + \tilde\phi(v) = a|v|^2\ .\label{zdok} \end{equation}

To use this, observe that for any dual pair of convex functions $\eta$ and 
$\eta^*$, Young's
inequality say that $\eta(v) + \eta^*(w) \ge v\cdot w$. Hence for all $v$ and 
$w$,
$${1\over 2}|v-w|^2 \ge {1\over 2}|v|^2 + {1\over 2}|w|^2 - \eta(v) - 
\eta^*(w)\ .$$
Now if $G$ is any element of ${\cal E}$, and $\gamma_0$ is the optimal coupling 
between $G$ and $F_0$,
we have
\begin{eqnarray}
W_2^2(G,F_0) &=& \int_{\R^d \times \R^d} {1\over 2}|v-w|^2\gamma_0({\rm d}v,{\rm 
d}w)\label{zzww} \\
&\ge& d\theta -  \int_{\R^d} \eta(v)G(v){\rm d}v  - \int_{\R^d} 
\eta^*(w)F_0(w){\rm 
d}w.\nonumber\end{eqnarray}
In the same way, we deduce that for any other dual pair of convex functions 
$\zeta$ and $\zeta^*$,
\begin{eqnarray}
W_2^2(G,F_1) &=& \int_{\R^d \times \R^d} {1\over 2}|v-w|^2\gamma_1({\rm d}v,{\rm 
d}w)\label{zzqq}\\
&\ge& d\theta -  \int_{\R^d} \zeta(v)G(v){\rm d}v  - \int_{\R^d} 
\zeta^*(w)F_1(w){\rm d}w\ 
.\nonumber  \end{eqnarray}
We now choose $\eta = \phi$ and $\zeta = \tilde\phi$. Then adding \pref{zzww} and 
\pref{zzqq}, and
on account of \pref{zdok},
\begin{eqnarray}
\quad\Phi(G) &=& W_2^2(G,F_0) + W_2^2(G,F_1)\label{rrzz}\\
&\ge &2d\theta - \int_{\R^d}\left(\phi(v) + \tilde \phi(v)\right)G(v){\rm d}v
\nonumber\\
&&- \int_{\R^d}\phi^*(w)F_0(w){\rm d}w - \int_{\R^d}\tilde \phi^*(w)F_1(v){\rm 
d}w\nonumber\\
&=& \left(2 - a\right)d\theta - \int_{\R^d}\phi^*(w)F_0(w){\rm d}w - 
\int_{\R^d}\tilde 
\phi^*(w)F_1(v){\rm d}w
\ .\nonumber \end{eqnarray}

Now suppose that $G=G_0$. Then for $\gamma_0$-almost every $(v,w)$, we have 
that  
$v\cdot w = \phi(v) + \phi^*(w)$ so that
$${1\over 2}|v-w|^2 = {1\over 2}|v|^2 + 
{1\over 2}|w|^2 - \phi(v) - \phi^*(w)
$$
and hence there is equality in \pref{zzww} when $G=G_0$ and $\eta = \phi$. In the 
same way, there is equality
in \pref{zzqq} when $G=G_0$ and $\zeta = \tilde\phi$.
Thus, the lower bound in \pref{rrzz} is saturated for $G = G_0$, and is in any 
case independent of $G$.
This proves that $G_0$ is the minimizer. 

It is now easy to compute the minimizing value. 
Theorem 3.1 tells us that $G_0(v) = a^dF_{1/2}(av)$ where $a$ depends only 
on 
$W_2^2(F_0,F_1)$, and is given explicitly by \pref{aform}.
Then, with this choice of $a$, 
$${1\over a}\nabla \psi_{1/2}\# F_0 = G_0\ .$$
Expressing this directly in terms of $\psi$ and computing in the familiar way, 
one finds
\begin{equation} W_2^2(F_0,G_0)  
= {d\theta\over a}\left[(a-1) + {W_2^2(F_0,F_1)\over 2d\theta} \right] = 
d\theta(1-a)\ .\hskip.5in
\label{vop} \end{equation}
Clearly, $W_2^2(F_0,G_0) =W_2^2(G_0,F_1)$, and so doubling the right-hand side 
of \pref{vop}
and inserting our formula for $a$, we obtain \pref{gonzo}.
Finally simple calculations confirm that $f$ is increasing and convex on 
$[0,1]$. \enddemo

We are now  prepared to consider discrete approximations to geodesics in~${\cal 
E}$. Let
${\cal G}$ be the set of continuous maps
$t\mapsto G_t$ from $[0,1]$ to ${\cal E}$ with $G_0 = F_0$ and $G_1 = F_1$. 

For each
natural number $k$, let ${\cal G}_k(F_0,F_1)$ denote the set of sequences
\begin{equation} \{G_0, G_1,\dots,G_{2^k}\}\label{gkdef} \end{equation}
where each $G_j$ is in ${\cal E}$,  $G_0 = F_0$, $G_{2^k} = F_1$,
and finally 
\begin{equation} W_2^2(G_{j+2},G_{j+1}) = W_2^2(G_{j+1},G_{j})\ \label{eqeq} \end{equation}
for all $j = 0, 1, \dots,2^k-2$.

For any path $t\mapsto G_t$ in ${\cal G}$ and any $k$, we obtain a sequence in 
${\cal G}_k(F_0,F_1)$ by an appropriate selection of times $t_j$ and by setting 
$G_j = G(t_j)$.

We next obtain a particular element 
$\{F^{(k)}_0,F^{(k)}_1,\dots,F^{(k)}_{2^k}\}$
of ${\cal G}_k(F_0,F_1)$ by successive midpoint projections
onto 
${\cal E}$ as follows: For $k=1$, let $F^{(1)}_0 = F_0$ and $F^{(1)}_2= F_1$ as 
we must. Define $F^{(1)}_1$
to be the midpoint of the chordal geodesic from $F_0$ to $F_1$, projected onto 
${\cal E}$ as in Theorem 3.3.
Then, supposing $\{F^{(k)}_0,F^{(k)}_1,\dots,F^{(k)}_{2^k}\}$ to be defined, put
$F^{(k+1)}_{2j} = F^{(k)}_j$ for $j = 0,1,\dots ,2^k$. Also, for $j = 0,1,\dots 
,2^k-1$,
let $F^{(k+1)}_{2j+1}$ be the midpoint of the chordal geodesic from 
$F^{(k)}_j$
to  $F^{(k)}_{j+1}$, projected onto ${\cal E}$ as in Theorem 3.3.

\proclaimtitle{Discrete geodesics} 
\proclaim{Lemma} For all $k\ge 1${\rm ,}
$$
\sum_{j=0}^{2^k-1}W_2(F^{(k)}_{j},F^{(k)}_{j+1}) \le 
\sum_{j=0}^{2^k-1}W_2(G_{j},G_{j+1})$$
for any $\{G_0, G_1,\dots,G_{2^k}\}$ in ${\cal G}_k(F_0,F_1)${\rm ,} and there is 
equality  when and only when
$$\{G_0, G_1,\dots,G_{2^k}\} = \{F^{(k)}_0,F^{(k)}_1,\dots,F^{(k)}_{2^k}\}\ .$$
\endproclaim

\demo{Proof} By   condition \pref{eqeq}, 
\begin{equation} \sum_{j=0}^{2^k-1}W_2(G_{j},G_{j+1}) = 
\left(\sum_{j=0}^{2^k-1}{W_2^2(G_{j},G_{j+1})
\over 2^{-k}}\right)^{1/2}\ .\label{soso} \end{equation}
We now claim that
$$ \sum_{j=0}^{2^k-1}W_2^2(F^{(k)}_{j},F^{(k)}_{j+1}) \le 
\sum_{j=0}^{2^k-1}W_2^2(G_{j},G_{j+1})$$
 and there is equality exactly when
$\{G_0, G_1,\dots,G_{2^k}\}$ = $\{F^{(k)}_0,F^{(k)}_1,\dots,F^{(k)}_{2^k}\}$. 
On 
account of \pref{soso},
once this is established,
the proof is complete.

For $k=1$, this is  implied by Theorem 3.3. 
For $k>1$, consider any\break $2^k+1$-tuple $\{G_0, G_1,\dots,G_{2^{k}}\}$ of 
elements of ${\cal E}$.
We {\it are not}\/ requiring\break
$\{G_0, G_1,\dots,G_{2^{k}}\}  \in {\cal G}_{k}$. The point is that we are going 
to reduce to the
case $k=0$ by successively erasing every other element. Even if  
$W_2(G_{j},G_{j+1}) = W_2(G_{j+1},G_{j+2})$
for all $j$, it is not necessarily the case that  $W_2^{\phantom{|}}(G_{j},G_{j+2}) = 
W_2^{\phantom{|}}(G_{j+2},G_{j+4})$ for all $j$,
so that the procedure of ``erasing midpoints'' does not take us from ${\cal G}_{k}$ 
to ${\cal G}_{k-1}$

Nonetheless, without assuming that $\{G_0, G_1,\dots,G_{2^{k}}\}\in {\cal 
G}_{k}$, we have from
Theorem 3.3, with $f$ given by \pref{gonzo}, \pagebreak that
\begin{eqnarray}
&&\label{pozq}\\
\sum_{j=0}^{2^{k}-1}W_2^2(G_j,G_{j+1}) &= &
\sum_{\ell=0}^{2^{k-1}-1}\left(W_2^2(G_{2\ell},G_{2\ell+1}) + 
W_2^2(G_{2\ell+1},G_{2\ell+2})\right)\nonumber\\
&\ge& \sum_{\ell=0}^{2^{k-1}-1} f(W_2^2(G_{2\ell},G_{2\ell+2}))\nonumber \\
&=&2^{k-1}\left({1\over 2^{k-1}}\sum_{\ell=0}^{2^{k-1}-1} 
f(W_2^2(G_{2\ell},G_{2\ell+2}))\right)\nonumber \\
&\ge&2^{k-1}f\left({1\over 2^{k-1}}\sum_{\ell=0}^{2^{k-1}-1} 
W_2^2(G_{2\ell},G_{2\ell+2})\right)
\nonumber \end{eqnarray}
where the last inequality is the convexity of $f$. 

Notice that both inequalities are saturated if and only if for each $\ell$, 
$G_{2\ell+1}$ is the
projected midpoint of the chordal geodesic connecting $G_{2\ell}$ and 
$G_{2\ell+2}$.

The proof is now easy to complete. Define a sequence $\{A_j\}$ inductively by 
$A_0 = W_2^2(F_0,F_1)$
and
\begin{equation} A_{j+1} = 2^jf\left(2^{-j}A_j\right)\ .\label{rolz} \end{equation}
Because these inequalities are  saturated for 
$\{G_0, G_1,\dots,G_{2^k}\}$ = $\{F^{(k)}_0,F^{(k)}_1,\break\dots,F^{(k)}_{2^k}\}$,
$$A_k = \sum_{j=0}^{2^{k}}W_2^2(F^{(k)}_j,F^{(k)}_{j+1})\ .$$
But a simple induction argument based on \pref{pozq} shows that 
$$\sum_{j=0}^{2^{k}}W_2^2(G_j,G_{j+1}) \ge A_k$$ with equality only in the stated 
case.
\enddemo

We can now define the distance ${\cal W}_2(F_0,F_1)$ on ${\cal E}$ induced by 
the $2$-Wasserstein metric:
\begin{equation} {\cal W}_2(F_0,F_1) = \lim_{k\to 
\infty}\sum_{j=1}^{2^k-1}W_2(F^{(k)}_{j},F^{(k)}_{j+1})\label{www} \end{equation}
where clearly the sequence on the right   in \pref{www} is increasing. In fact, 
Lemma 3.4 tells us that the geodesic from
$F_0$ to $F_1$ on ${\cal E}$ is obtained by the following simple rule: Take the 
chordal geodesic 
$t\mapsto F_t$ from
$F_0$ to $F_1$ in ${\cal P}$, and rescale each $F_t$ onto ${\cal E}$ as in 
Theorem 3.1. Then reparametrize
this path in ${\cal E}$ so that it runs at constant speed. This is the geodesic. 
Note that this same procedure
produces geodesics on the sphere $S^{d-1}$ in $\R^d$.

It is now an easy matter to compute the distance ${\cal W}_2(F_0,F_1)$. One way 
is to compute $\lim_{k\to\infty}A_k$
for the sequence given by $A_0 = W_2^2(F_0,F_1)$ and \pref{rolz}. This is  
straightforward; it is easy to recognize the
iteration as the same iteration one gets by dyadically rectifying an arc of the 
circle.

We find it more enlightening  to obtain an explicit parametrization
of the corresponding geodesic, and to use the Riemannian metric for the\break
$2$-Wasserstein distance.

To begin the computation, let $\psi$ be the convex function such that 
$\nabla\psi\#F_0\break = F_1$.
We may assume without loss of generality that $u=0$; this will simplify the 
computation.
Then define $F_t$ as in \pref{geoA1} and \pref{geogeo}, and let $\tilde F_t$ be 
the projection of $F_t$ onto ${\cal E}$
as in Theorem 3.1. Since $u=0$,
$$\tilde F_t = \left({1\over a(t)}\nabla \psi_t\right)\#F_0$$
where $\psi_t$ is defined in terms of $\psi$ as usual
and where
$$a(t) = \sqrt{1 -4t(1-t){W_2^2(F_0,F_1)\over 2d\theta}}\ .$$
Notice that the gradient vector field on $\R^d$ that represents the tangent 
vector $\partial \tilde F_t/
\partial t$ has two terms: One is a rescaling of the gradient vector field on 
$\R^d$ that represents
$\partial F_t/ \partial t$, and the other generates a dilation to keep the path 
on~${\cal E}$. 

Next, we have from Theorem 2.3 that for any test function $\chi$ on $\R^d$, 
after some computation,
\begin{eqnarray*}
{{\rm d}\over {\rm d}t}\int_{\R^d}\chi(v)\tilde F_t(v){\rm d}v &= &
{{\rm d}\over {\rm d}t}\int_{\R^d}\chi\left({v\over a(t)}\right)F_t(v){\rm 
d}v \\ 
&=&\int_{\R^d}\nabla \chi(v)\cdot \left({1\over a(t)}\nabla\eta_t(a(t)v)
- {\dot a(t)\over a(t)}v\right)\tilde  F_t(v){\rm d}v\ ,
\end{eqnarray*}
where $\eta_t$ is given by \pref{etadef3}. Hence, from \pref{metric}, we have 
\begin{eqnarray*}
g\left({\partial \tilde F_t\over \partial t},{\partial \tilde F_t \over \partial 
t}\right) 
&=&{1\over 2}\int_{\R^d}\left|{1\over a(t)}\nabla\eta_t(a(t)v)
- {\dot a(t)\over a(t)}v \right|^2\tilde F_t(v){\rm d}v \\
&=&{1\over 2a^2(t)}\int_{\R^d}\left|\nabla\eta_t(v)
- {\dot a(t)\over a(t)}v \right|^2 F_t(v){\rm d}v\ .\end{eqnarray*}
By \pref{rice6},
$\int_{\R^d}\left|\nabla\eta_t(v)\right|^2 F_t(v){\rm d}v   =2 W_2^2(F_0,F_1)$,
and clearly
$\int_{\R^d}\left|v\right|^2 F_t(v){\rm d}v = a^2(t)d\theta$.
Finally, by Theorem 2.3 and familiar computations,
\begin{eqnarray*}
&&
\int_{\R^d}\left(\nabla\eta_t(v)\cdot v\right)F_t(v){\rm d}v\\
 &&\qquad =
{1\over 2t}\int_{\R^d}\left(   |\nabla(\psi_t)^*(v) -v|^2 + |v|^2 - 
|\nabla(\psi_t)^*(v)|^2 \right)F_t(v){\rm d}v \\
 &&\qquad ={1\over 2t}\left(2W_2^2(F_0,F_t) + (a^2(t) -1)d\theta\right) = (2t-1) 
W_2^2(F_0,F_1)\ .\pagebreak \end{eqnarray*}
Putting all of this together,   one has,
after some algebra,
\begin{eqnarray*}
g\left({\partial \tilde F_t\over \partial t},{\partial \tilde F_t \over \partial 
t}\right) 
&=&
{1\over 2a^2(t)}\left[2W_2^2(F_0,F_1)  + \left({\dot a(t)\over 
a(t)}\right)^2a^2(t)d\theta\right.
\\
&&  \left. \phantom{\left({\dot a(t)\over 
a(t)}\right)^2} -2
{\dot a(t)\over a(t)}(2t-1)W_2^2(F_0,F_1) \right] \\
&=&W_2^2(F_0,F_1) {1\over a^4(t)}\left[1 - {W_2^2(F_0,F_1)\over 2d\theta}\right]\ 
.\end{eqnarray*}

Now we reparametrize to achieve constant unit speed. We take the map
$t\mapsto \tau(t)$ to be differentiable and increasing. Then with $\tilde F_\tau 
=  
\tilde F_{\tau(t)}$,
\begin{equation} 1 = g\left({\partial\tilde F_\tau\over \partial \tau},
{\partial\tilde F_\tau\over \partial \tau}\right) =
g\left({\partial\tilde F_t\over \partial t},
{\partial\tilde F_t\over 
\partial t}\right)\left|{\rm d}t\over {\rm d}\tau\right|^2\label{normalize} \end{equation}
provided
$${{\rm d}\tau(t)\over {\rm d}t} = W_2(F_0,F_1) {1\over a^2(t)}\sqrt{1 - 
{W_2^2(F_0,F_1)\over 2d\theta}}
\ .$$
This is solved by
$$\tau(t) = \sqrt{d\theta\over 2}{\rm 
arctan}\left((2t-1)\sqrt{{W_2^2(F_0,F_1)\over
2d\theta -W_2^2(F_0,F_1)}}\right)
$$ for which $\tau(1/2) = 0$ and
\begin{equation} {\cal W}_2(F_0,F_1) =\tau(1) - \tau(0)=  2\sqrt{d\theta\over 2}{\rm 
arctan}\left(\sqrt{{W_2^2(F_0,F_1)\over
2d\theta -W_2^2(F_0,F_1)}}\right).\hskip.25in\label{dist} \end{equation}
This has a very simple interpretation: Consider two points on a circle of radius 
$R$, and let $D$ be
the length of the chord that they terminate. The arc joining them subtends an 
angle $2\phi$ where
$$\tan(\phi) = \sqrt{D^2\over 4R^2 - D^2}\ ,$$
and hence the length of the arc joining them is
\begin{equation} 2R{\rm arctan}\left(\sqrt{D^2\over 4R^2 - D^2}\right)\ .\label{xlan} \end{equation}
Since $\sqrt{(d\theta)/2}$ is the radius $R_\theta$ of ${\cal E}$,  in that this 
is the $2$-Wasserstein distance from any
point in ${\cal E}$ to the unit mass at $u$, and since $W_2(F_0,F_1)$ is the 
chordal separation of
$F_0$ from $F_1$ in the $2$-Wasserstein distance, we have that \pref{xlan}, with 
$R = \sqrt{(d\theta)/2}$
and $D = W_2(F_0,F_1)$, gives us ${\cal W}_2(F_0,F_1)$. It is somewhat simpler 
to express this in
terms of sines instead of tangents. From \pref{xlan} it is easy to deduce that
\begin{equation} W_2(F_0,F_1) = 
2R_\theta\sin\left({{\cal W}_2(F_0,F_1)\over 2R_\theta}\right)\ ,\label{geodist1} \end{equation}
\begin{equation} {\cal W}_2(F_0,F_1) = 
2R_\theta{\rm arcsin}\left({W_2(F_0,F_1)\over 2R_\theta}\right)\ 
.\label{geodist2} \end{equation}
We summarize this in the following theorem:

 \proclaimtitle{Geometry of ${\cal E}$}
\proclaim{Theorem}
 Let ${\cal W}_2(F_0,F_1)$ denote the distance between any two points $F_0$ 
and $F_1$ of ${\cal E}$
in the metric induced on ${\cal E}$ by the\break $2$\/{\rm -}\/Wasserstein metric. Then ${\cal 
W}_2(F_0,F_1)$ is
related to $W_2(F_0,F_1)$ through {\rm \pref{geodist1}} and {\rm \pref{geodist2}.} Moreover{\rm ,} 
the geodesic on ${\cal E}$
between $F_0$ and  $F_1$ is obtained from the chordal geodesic  in ${\cal P}$ 
between $F_0$ and  $F_1$
by the following procedure\/{\rm :}\/ Let $t\mapsto F_t${\rm ,} $t\in [0,1]${\rm ,} denote the chordal 
geodesic. Then{\rm ,} for each
such $t${\rm ,} let $\tilde F_t$ denote the unique point in ${\cal E}$ that is closest 
to $F_t${\rm ,} which is simply obtained
from $F_t$ by dilating about the mean $u$. This path{\rm ,} reparametrized to run 
at constant speed{\rm ,} is the geodesic on 
${\cal E}$ between $F_0$ and  $F_1$.
\endproclaim

This theorem strongly encourages one to think of ${\cal E}$ in spherical terms, 
though we see from
\pref{lowvarbnd} that the chordal distance between any two points on ${\cal E}$
is no more than $\sqrt{2}$ times the radius of ${\cal E}$, as given by 
\pref{UU45}, as on the spherical cap
with the azimuthal angle $\phi$ ranging over $0\le\phi\le \pi/4$.

We apply this to deduce a criterion for displacement convexity on the 
constrained manifold ${\cal E}$.
We say that a functional  $\Phi$ is displacement convex on ${\cal E}$ in case 
for all geodesics
$t\mapsto G_t$ in ${\cal E}$, the function $t\mapsto \Phi(G_t)$ is convex. If 
the gradient vector field
$\nabla \eta$ on $\R^d$ is the tangent vector at $t=0$ to a geodesic 
$t\mapsto G_t$ in ${\cal E}$, we define
\begin{equation} {\cal H}{\it ess}\, \Phi(G_0)(\nabla \eta,\nabla \eta) = 
{{\rm d}^2 \over {\rm d}t^2}\Phi(G_t)\bigg |_{t=0}\ .\label{conhess} \end{equation}
This should be compared with \pref{hess}. The differences lie in the different 
classes of geodesics
being considered in the two cases, as well as the fact that 
\begin{equation} \int_{\R^d}v\cdot \nabla \eta(v)G_0(v){\rm d}v = 0\qquad{\rm and}\qquad 
\int_{\R^d}\nabla \eta(v)G_0(v){\rm d}v = 0\hskip.45in \label{musthave} \end{equation}
must hold for $\nabla \eta$ to represent a tangent vector to ${\cal E}$ at 
$G_0$.

Since we have determined the geodesics in ${\cal E}$, it is now a simple matter 
to determine a criterion 
for displacement convexity in ${\cal E}$.
\proclaimtitle{Displacement convexity in ${\cal E}$} 
\proclaim{Theorem}  Let $G\mapsto \Phi(G)$ be any 
functional of the form
$$\Phi(G) = \int_{\R^d} g(G(v)){\rm d}v$$ 
where $g$  is 
twice continuously differentiable on $\R_+$. Define the function $h$ 
by $h(t) = tg'(t) - g(t).$ Suppose that $F\in  {\cal E}_{u,\theta}$ is such that 
$h(F)$ is integrable{\rm ,} and that at $F${\rm ,} 
$$ G \mapsto {\rm Hess}\, \Phi(G)(\nabla\eta,\nabla\eta)
$$
is continuous in the $2$\/{\rm -}\/Wasserstein metric for all test functions $\eta$. Then 
\begin{eqnarray}
{\cal H}{\it ess}\, \Phi(F)(\nabla \eta,\nabla \eta) &= &
{\rm Hess}\, \Phi(F)(\nabla \eta,\nabla \eta)\label{supcon}\\
&&+ {d\over 2 R^2_\theta}
\left(\int_{\R^d}h(F){\rm d}v\right)\int_{\R^d}
|\nabla \eta|^2F{\rm d}v\ ,\nonumber \end{eqnarray}
where $R_\theta = \sqrt{d\theta/2}$ is the radius of ${\cal E}_{u,\theta}${\rm ,}
and $\nabla \eta$ is any  gradient vectorfield satisfying {\rm \pref{musthave}}
In particular{\rm ,} if $\Phi(F) = S(F)$ is the entropy $\int_{\R^d}\ln(F(v))F(v){\rm 
d}v$
of $F${\rm ,}
\begin{equation}
{\cal H}{\it ess}\, S(F)(\nabla \eta,\nabla \eta) = 
{\rm Hess}\,  S(F)(\nabla \eta,\nabla \eta) 
+ {d\over 2 R^2_\theta}\int_{\R^d}
|\nabla \eta|^2F{\rm d}v\ ,\hskip.25in\label{supentr} \end{equation}
and thus the entropy is uniformly  convex on the constrained manifold ${\cal 
E}_{u,\theta}$.
\enddemo

\demo{Proof} Without loss of generality, suppose $u=0$.
For any $F\in {\cal E}$, let $t\mapsto \tilde G_t$ be a
geodesic in ${\cal E}$ passing through $F$ with unit speed at $t=0$. Pick 
$\delta>0$  sufficiently small that $\tilde G_\delta$ and $\tilde 
G_{-\delta}$ are both defined.
By definition ${\cal W}_2^2(\tilde G_{-\delta},\tilde G_{\delta})= 4\delta^2$. 
Define $h>0$ by
$W_2^2(\tilde G_{-\delta},\tilde G_{\delta})= 4\ h^2$.
By Theorem 3.5,
\begin{equation} h = R_{\theta}\sin\left(\delta\over  R_{\theta}\right) = \delta + {\cal 
O}(\delta^3)\ .\label{8OP8} \end{equation}

Now let $t\mapsto G_t$ be the {\it chordal} geodesic, in ${\cal P}$, from 
$\tilde G_{-\delta}$ to $\tilde G_{\delta}$ parametrized so that
$\tilde G_{-\delta} = G_{-h}$ and $\tilde G_{\delta} = G_h$. By Theorem 3.3, 
$\tilde G_0 = F$ is obtained from $G_0$ by dilation:
\begin{equation} \tilde G_0(v) = a^dG_0(av)\label{MNQ3} \end{equation}
where
\begin{equation} a = \sqrt{1 - {h^2\over R_\theta^2}}\ .\label{9op9} \end{equation}
Now
\begin{eqnarray}  
&\hskip-10pt\hskip-10pt& \label{R9Z7}\\
{1\over \delta^2}\left[{1\over 2}\left(\Phi(\tilde G_{\delta})\! +\!   \Phi(\tilde 
G_{-\delta})\right) \! -  \! 
\Phi(\tilde G_{0})\right]
&\hskip-10pt=\hskip-10pt& {\Phi(G_{0})\!  -\!  \Phi(\tilde G_{0})\over \delta^2}\nonumber\\
&\hskip-10pt\hskip-10pt & - {1\over h^2}\left[{1\over 2}\left(\Phi( G_{h}) \! + \!  \Phi( G_{-h})\right) \! -\! 
\Phi( G_{0})\right]{h^2\over \delta^2}.\nonumber\end{eqnarray}
Next, since
$\Phi(G_{0}) - \Phi(\tilde G_{0}) = a^d\int_{\R^d}g(a^{-d}F(v)){\rm d}v - 
\int_{\R^d}g(F(v)){\rm d}v$, 
it follows from \pref{9op9} and the definition of $h$ that
\begin{equation} \lim_{\delta\to 0}{\Phi(G_{0}) - \Phi(\tilde G_{0})\over \delta^2} = {d\over 
2R^2_\theta}
\int_{\R^d} h(F(v)){\rm d}v\ .\label{R5y7} \end{equation}

By \pref{8OP8}, the continuity of ${\rm Hess}\, \Phi$ at $F$ and our previous 
definitions,
$$\lim_{\delta\to 0}{1\over h^2}\left[{1\over 2}\left(\Phi( G_{h}) +  \Phi( 
G_{-h})\right) -
\Phi( G_{0})\right]{h^2\over \delta^2} = {\rm Hess}\, \Phi(F)\ .$$
Combining this, \pref{R5y7} and \pref{R9Z7}, we obtain \pref{supcon} from which the 
rest of the
result easily follows. \enddemo

As an application, we deduce a strengthened form of an 
inequality due to Talagrand \cite{Tal}. Let $G_0$ be a Gaussian density in 
${\cal E}_{ \theta,u}$. Let $F$ be any other density in 
${\cal E}_{ \theta,u}$. Let $F_s$ be the geodesic in  ${\cal E}_{ \theta,u}$,
parametrized by arclength, starting at $F$ and going to $G_0$. Then by
\pref{supentr},
\begin{eqnarray*}
S(F)-S(G_0) &\hskip-8pt=\hskip-8pt& \int_0^{{\cal W}_2(G_0,F)}S'(F_s){\rm d}s \\
&\hskip-8pt=\hskip-8pt& \int_0^{{\cal W}_2(G_0,F)}\left(S'(G_0) + 
\int_0^s S''(F_r){\rm d}r\right){\rm d}s
\ge {1\over 2} {d\over 2 R^2_\theta}{\cal W}_2^2(G_0,F)\ .\end{eqnarray*}
We have used the fact that $S'(G_0)= 0$ since $S(F) \ge S(G_0)$ by 
the entropy-minimizing property of Gaussians. Also, since both $F$ and $G_0$
lie in ${\cal E}_{ \theta,u}$,\break
$S(F) -S(G_0)= H(F|G_0)$, the relative entropy of $F$ with respect to $G_0$.
Therefore, since $R^2 = 2/(d\theta)$,
$$H(F|G_0) \ge 
{1\over 2\theta}{\cal W}_2^2(G_0,F)\ ,$$
which is Talagrand's inequality, except that here ${\cal W}_2^2(G_0,F)$ replaces 
the smaller quantity $W_2^2(G_0,F)$.

\section{The Euler-Lagrange equation} 

For fixed $h>0$, and a given density $F_0\in {\cal E}_{\theta,u}$, 
we seek to minimize
the functional
\begin{equation}
I(F) = \left[{W_2^2(  F_0,F)
\over \theta} +hS(F)\right]\ \label{fudef}, \end{equation}
subject to the constraint that $F\in {\cal E}_{\theta,u}$.
 
This functional is
strictly convex and our constraints are convex, and hence if any minimizer does 
exist, it would also be unique. The
existence issue will be settled in the next section. Here we shall
derive the Euler Lagrange equation that would be satisfied by any  minimizer
in our variational problem, and derive some consequences of satisfying this 
equation.

\proclaim{Theorem} Suppose that $F_1$ is a minimizer of
the functional given in {\rm \pref{fudef}} subject to the constraint that $F_1$
has the same mean and variance as 
$F_0$. Let $ \psi$ be the convex function on $\R^d$ 
such that 
\begin{equation} \nabla   \psi\# F_1 =   F_0\ .\label{pel} \end{equation}
Then
\begin{equation} \int_{\R^d}|\nabla \ln F_1|^2F_1(v){\rm d}v < \infty\ \label{nel} \end{equation} 
and
\begin{equation} \nabla   \psi(v) = v + h\theta\nabla_v\left(\ln {F_1\over M_{F_1}}\right)
+(u - v)\left[{W^2_2(F_1,  F_0)\over d\theta} \right] 
\label{mel} \end{equation}
where
for any $F\in {\cal P}${\rm ,} $M_F$ denotes the isotropic Gaussian density with the 
same mean and variance as $F$.
\endproclaim
 
\demo{Proof} Consider a function $\xi: \R^d\rightarrow \R^d$ satisfying
\begin{equation} \int_{\R^d}\xi(v)F_1(v){\rm d}v  =0\qquad{\rm and}\qquad
\int_{\R^d}(\xi(v)\cdot v)F_1(v){\rm d}v  =0\ . \label{chicond} \end{equation}
Then  define the flow
$T_t(v) = v + t\xi(v)$ and the curve of densities
$G(t) =T_t\# F_1$. Finally, let $\tilde G(t)$ be the projection of $G(t)$ onto
${\cal E}$ as in Theorem 3.1. 
Let $u_1$ and $d\theta_1$ be the mean and variance of $F_1$. Then
by Theorem 3.1,
$\tilde G(t,v) = a(t)^d G(t, a(t)(v-u(t)) + u_1)$,
where, by \pref{chicond}
\begin{equation} a(t) = 1 + {\cal O}(t^2) \qquad{\rm and}\qquad  u(t) = u_1+ {\cal O}(t^2)\ 
.\label{NMD2} \end{equation} 
We can also write $\tilde G(t) = \tilde T_t\# F_1$ where
$\tilde T_t(v) = \left(v + t\xi\left(v/a(t)\right)\right/a(t)$.

The argument here is adapted from the corresponding argument in \cite{JKO}. 
First, consider the entropy. By direct calculation and \pref{NMD2},
$$S(\tilde G(t)) - S(\tilde G(0))= - t\int_{\R^d} F_1(v) \nabla \cdot \xi(v) 
{\rm 
d}v + {\cal O}(t^2)$$
and so
$$\lim_{t\to 0^+}{S(\tilde G(t))-S(F_1)\over t} = -\int_{\R^d} 
 F_1(v)\nabla \cdot \xi(v){\rm d}v\ .$$

To compute the variation in the $2$-Wasserstein distance, note that
since $\tilde T_t\# F_1 = \tilde G(t)$, $\nabla   \psi \circ \tilde T_t^{-1}\# 
\tilde G(t) =  F_0$. Thus
\begin{eqnarray*}
W_2^2(\tilde G(t),  F_0) &\le &{1\over 2}\int_{\R^d}
|\nabla   \psi \circ \tilde T_t^{-1}(v) - v|^2 \tilde G(t,v){\rm d}v \\
&=& {1\over 2}\int_{\R^d}
|\nabla   \psi  - \tilde T_t(v)|^2 F_1(v){\rm d}v \\
&\le& W_2^2(F_1,  F_0) - 
t\int_{\R^d}\left(\nabla   \psi - v\right)\cdot \xi F_1(v){\rm d}v + o(t)
\ .\end{eqnarray*}
Now it follows easily  that
\begin{equation} \limsup_{t\to 0^+}{W_2^2(G(t),F_0)-W_2^2(F_1,  F_0)\over t} \le \int_{\R^d} 
\left( v- \nabla   \psi(v) \right)F_{1}(v)\cdot \xi(v){\rm d}v\ .\label{sec2B} \end{equation}

We deduce that
$$
\int_{\R^d}\left( 
\left(\nabla   \psi(v)-v \right){F_1(v)\over \theta}
\right)\cdot \xi(v){\rm d}v \le  -h\int_{\R^d}  F_1(v) \nabla \cdot\xi(v)
$$
for all smooth and compactly supported $\xi$ satisfying \pref{chicond}. Since 
these conditions
are still satisfied if $\xi$ is replaced by $-\xi$, we have that
$$
\int_{\R^d}\left( 
\left(\nabla   \psi(v)-v \right){F_1(v)\over \theta}
\right)\cdot \xi(v){\rm d}v = - h\int_{\R^d}  F_1(v) \nabla \cdot\xi(v)
$$
for all smooth and compactly supported $\xi$ satisfying \pref{chicond}.
Hence
\begin{equation} \left( 
\left(\nabla   \psi(v)-v \right){F_1(v)\over \theta}
- h\nabla F_1(v)\right) = (A + B(u - v))F_1(v)\label{elone} \end{equation}
for some vector $A$ and scalar $B$. It follows from this that \pref{nel} holds.

Integrating both sides of \pref{elone} in $v$, one learns that $A=0$.
If one takes the inner product of both sides with $(u - v)$,
and then integrates, one learns
$d\theta B =  W_2^2(F_1,  F_0)/\theta  - dh$
since
$$\int_{\R^d}\left(\nabla   \psi(v)-v \right)\cdot v F_1(v){\rm d}v = 
W_2^2(F_1,  F_0)\ .$$
Combining this and \pref{elone}, we obtain \pref{mel}. \enddemo

Now still assuming that the minimizer $F_1$ exists, we ask {\it what properties 
does $F_1$ inherit
from} $F_0$? We shall show, using the fact that $F_1$ satisfies the 
Euler-Lagrange equation \pref{mel}
and \pref{pel}, that $F_1$ inherits some localization properties from $F_0$.
Specifically, let $\zeta$ be a nonnegative, increasing convex function on 
$\R_+$ with the property that 
$\lim_{t\to\infty}\zeta(t)/t = \infty$ and that $\zeta(0)=0.$ Suppose that
\begin{equation} \int_{\R^d}\zeta(|v|^2)F_0(v){\rm d}v = C < \infty\ .\label{4tu} \end{equation}
This quantity provides a quantitative measure of the localization of 
$|v|^2F_0(v)$ in that
$$\int_{|v|^2>t}|v|^2F_0(v){\rm d}v \le {t\over\zeta(t)}C\ ,$$
and the right-hand side tends to zero as $t$ increases. Here, we have used that 
$t \rightarrow \zeta(t)/t$ is nondecreasing. If we knew that $F_1$ satisfied the 
same inequality,
we would have a quantitative localization estimate on $F_1$. We shall see below 
that this is almost the case:
The function $\zeta$ is modified slightly in passing from $F_0$ to $F_1$.

First, we need to explain where the original $\zeta$ comes from. We could take 
$\zeta(t) = (1+t)^{1+\varepsilon}$
if we assumed that $F_0$ possessed more than second moments. Since we wish to 
make a statement about
generic elements $F_0$ of ${\cal E}_{u,\theta}$, we use  a minor variant of 
a lemma of de~la~Vall\'ee-Poussin, which says that
for any probability density
$F_0$ with $\int_{\R^d}|v|^2F_0(v){\rm d}v < \infty$, there is a
a nonnegative, increasing convex function on $\R_+$ with the property that 
$\lim_{t\to\infty}\zeta(t)/t = \infty$ such that \pref{4tu} holds, and finally, 
that $\|\zeta''\|_\infty \le 1$. Everything up to the last condition is 
standard, 
though the usual construction of $\zeta$ is such that $\zeta''$ is a series of 
Dirac masses.
We therefore sketch a short proof. Without loss of generality, we may suppose 
that $u=0$ and $\theta=1/d$.

Let 
$$\lambda(t) = \int_{|v|^2 >t}F_0(v){\rm d}v \qquad{\rm and}\qquad 
\mu(t) = \int_{|v|^2 >t}|v|^2F_0(v){\rm d}v $$
so that
$1 = \int_{\R^d}|v|^2F_0(v){\rm d}v =  \int_0^\infty \lambda(t){\rm d}t$, and 
that
\begin{equation} \mu(t) = \int_t^\infty \lambda(u){\rm d}u + \lambda(t) \ge 
\sum_{n>t}\lambda(n)\ 
.\label{sel} \end{equation} 
Here, we have used the layer cake representation theorem. Now define $t_k$ by 
$t_0 =0 $ and for $k\ge 1$, $t_k = \inf\{ t \ |\ \mu(t) < 2^{-k}\}$.  
Since $F_0(v){\rm d}v$ is absolutely continuous, $\mu(t_k) = 2^{-k}$. Then by 
\pref{sel},
\begin{equation} 1 = \sum_{k=1}^\infty \mu(t_k) = \sum_{k=1}^\infty \sum_{n>t_k}\lambda(n) = 
\sum_{n=1}^\infty
g(n)\lambda(n)\label{tel} \end{equation}
where $g(0)=0$ and for all $n \ge 1$, $g(n) = \max\{k \ | \ t_k < n\}$.
Clearly, $\lim_{n\to\infty}g(n) =\infty$ and $g(n+1) \ge g(n).$  
Next, set $h(0)=0$  and for $n \ge 1,$ define $h(n)$ recursively by $h(n) - 
h(n-1) =1$ if $g(n) - g(n-1)>0$, and 
$h(n) - h(n-1) =0$ otherwise. Then 
$$ h(n) = \sum_{k=1}^n(h(k) - h(k-1)) \le \sum_{k=1}^n(g(k) - g(k-1)) = g(n)
$$
but also clearly $\lim_{n\to\infty}h(n) =\infty$ since $g(n)$ must increase 
infinitely often.

Now define $h(t)$ for all $t>0$ by linear 
interpolation of $h(n)$, and then define 
$\zeta(t) = \int_0^th(s){\rm d}s$.
Note that $\zeta(t)$ is a
continuously differentiable convex increasing function with $\|\zeta''\|_\infty 
\le 1$, and $\lim_{t\to\infty}\zeta(t)/t = \infty$.
Also, since $\zeta(t)$ is increasing and $\lambda(t)$ 
is decreasing,
$$\int_0^\infty \zeta'(t)\lambda(t){\rm d}t \le \sum_{n=0}^\infty
h(n+1)\lambda(n) \le \sum_{n=0}^\infty
(1+g(n))\lambda(n) \le 3\ ,$$ where the last inequality follows from \pref{tel}.
Since
$\int_{\R^d}\zeta(|v|^2)F_0(v){\rm d}v =  \int_0^\infty \zeta'(t) 
\lambda(t){\rm d}t $,
\pref{4tu} holds.

We are now ready to prove the following:

\proclaim{Theorem} Suppose $F_0$ is any element of ${\cal 
E}_{\theta,u}${\rm ,} 
and suppose $\psi$ is a convex 
potential with $\nabla\psi\# F_1 = F_0$ such that $\psi$ and $F_1$ satisfy 
{\rm \pref{mel}. }
Then there are a nonnegative{\rm ,} increasing convex  function $\zeta(t)$ such that 
$\lim_{t\to\infty}\zeta(t)/t= \infty$ and $\|\zeta''\|_\infty \le 1${\rm ,} and a 
finite constant $C${\rm ,} both depending only on $F_0${\rm ,} so that
$$\int_{\R^d}\zeta(\alpha |w-u|^2)F_1(w){\rm d}w < C$$
for some $\alpha$ depending only on $h,$ $W_2(F_0,F_1),$ and $\theta$.
\endproclaim

\demo{Proof}
Without loss of generality, we continue to assume that $u=0$ and $\theta=1$, and 
thus
$$ \nabla\psi(w) = \alpha w +h\nabla\ln(F_1(w))
$$
for some constant $\alpha>0$ that is readily computed from \pref{mel}. Now let 
$\zeta(t)$ be the increasing  convex function provided by the 
variant of the de~la~Vall\'ee-Poussin lemma. Then, $v \rightarrow \zeta(|v|)$ 
is convex and so, 
\begin{eqnarray*}
\int_{\R^d}\zeta(|v|^2)F_0(v){\rm d}v &\hskip-8pt =\hskip-8pt &
\int_{\R^d}\zeta(|\nabla\psi(w)|^2)F_1(w){\rm d}w \\
&\hskip-8pt=\hskip-8pt&\int_{\R^d}\zeta(|\alpha w +h\nabla\ln(F_1(w))|^2)F_1(w){\rm d}w \\
&\hskip-8pt\ge\hskip-8pt&\int_{\R^d}\zeta(\alpha|w|^2)F_1(w){\rm d}w +
2h\alpha\int_{\R^d}\zeta'(\alpha|w|^2)w\cdot \nabla F_1(w){\rm d}w \\
&\hskip-8pt=\hskip-8pt&\int_{\R^d}\zeta(\alpha|w|^2)F_1(w){\rm d}w -
2h\alpha^2\int_{\R^d}\zeta''(\alpha|w|^2)|w|^2\cdot  F_1(w){\rm d}w 
 \\
&\hskip-8pt\hskip-8pt&- 2h\alpha d\int_{\R^d}\zeta'(\alpha|w|^2)\cdot  F_1(w){\rm d}w\ .
\end{eqnarray*}
Since $\int_{\R^d}|w|^2\cdot  F_1(w){\rm d}w=1$,
\begin{equation} \int_{\R^d}\zeta(\alpha|w|^2)F_1(w){\rm d}w \le 
\int_{\R^d}\zeta(|v|^2)F_0(v){\rm d}v
+ 2h\alpha^2\left(1 + d\right)\ ,\hskip.25in \label{4wingo} \end{equation}
where we are using the fact that $\|\zeta''\|_\infty \le 1$ and $\zeta'(t) \le 
t$ when
$\zeta$  is the function provided by the above variant of the 
de~la~Vall\'ee-Poussin
lemma. \enddemo

\section{Existence of minimizers}

To simplify the notation, we fix $u=0$ and $\theta =1$ throughout this 
section. The main goal  is to prove that a minimizer exists for 
\pref{fudef}. As explained in the introduction, it suffices to find a density
$F_1\in {\cal E}$ and a convex potential $\psi$ with $\nabla\psi\#F_1= F_0$ such 
that the Euler-Lagrange equation \pref{mel} is satisfied.

In this, we make essential use of the dual version of the variational
characterization of the $2$-Wasserstein metric. This says that for all $F_0$ and 
$F$ in~${\cal E}$,
\begin{eqnarray}
&& \label{dual}\\
 d - W_2^2(F_0,F)
& \hskip-8pt=\hskip-8pt&\inf\left\{ \int_{\R^d} \phi(v)F_0(v){\rm d}v\right. \nonumber\\
&\hskip-8pt\hskip-8pt&\qquad \left.+ \int_{\R^d} \psi(w)F(w){\rm 
d}w\ \bigg|\ 
\phi(v)+\psi(w) \ge v\cdot w \ {\rm a.e.}\ \right\}\ ,\nonumber \end{eqnarray}
where `almost everywhere' refers to the measure $F_0(v)F_1(w){\rm d}v{\rm d}w$.
Furthermore, the minimizing pair, which exists, consists of a dual pair of 
convex functions. 
That is, 
we may assume that $\phi$ and $\psi$ are Legendre transforms of one another. The 
gradients of 
the minimizing pair
provide the optimal transport plans; i.e., $\nabla \phi\# F_0 = F$ and
$\nabla \psi\# F = F_0$. A good reference for this is \cite{Bren} or 
\cite{Evans}.

We shall assume strong assumptions on $F_0\in {\cal E}$, 
which we shall later remove; namely 
we suppose that $F_0$ is supported in $B_R$, the centered ball of radius $R$, 
and that on $B_R$ it is bounded below by some
strictly positive number $\alpha$. Then for any other density $F$ in ${\cal P}$, 
these hypotheses
impose some regularity on the optimal map $\nabla \psi \#F = F_0$. In 
particular, 
\begin{equation} |\nabla \psi(v)| \le R\label{5lip} \end{equation}
for all $v$, which means that $\psi$ is Lipschitz.

Now define $\eta(t)$ by  
$$
\eta(t)= \left\{  \begin{array}{ll} +\infty & \hbox{if  } t< 0 , \\  t\ln t & \hbox{if } t\ge0. \end{array}
\right.
$$
Then the Legendre transform $\eta^*(s)$ of $\eta(t)$ is
$\eta^*(s) = e^{s-1}$.
We shall use use the notation $\eta^*$ throughout this section to emphasize the 
fact that we do not make much use of
the specific form of $\eta$ in our analysis; this point is discussed further at 
the end of the section.
Then
$$S(F) = \int_{\R^d} \eta(F){\rm d}v\ ,$$
and for any dual convex pair of functions $\phi$ and $\psi$,
\begin{equation}  I(F) \ge hS(F) +d - \left(\int_{\R^d} \phi(v)F_0(v){\rm d}v + \int_{\R^d} 
\psi(w)F(w){\rm d}w\right)\ ,\hskip.5in\label{5upper} \end{equation}
where $I(F)$ is given by \pref{fudef}. Moreover, by Young's inequality, $\eta(t) 
+ \eta^*(s) \ge st$, and thus we have that for any
$a\in R^d$ and any $b\in \R$,
\begin{equation} \eta(F) + \eta^*\left({a\cdot w + b|w|^2/2 + \psi(w)\over h}\right)\ge
{a\cdot w + b|w|^2/2 + \psi(w)\over h}F\ .\hskip.25in\label{5young} \end{equation}
Integrating yields
\begin{equation} hS(F) - \int_{\R^d}\psi(w)F(w){\rm d}w \ge  {d \over 2}b - h\int_{\R^d} \eta^*
\left({a\cdot w + b|w|^2/2 + \psi(w)\over h}\right)
{\rm d}w\ . \label{5aaa} \end{equation}
Therefore, introduce the functional
\begin{equation} J(a,b,\phi,\psi) =   d -\int_{\R^d} \phi(v)F_0(v){\rm d}v + {d \over 2}b 
- h\int_{\R^d} \eta^*\left({a\cdot w + b|w|^2/2 + \psi(w)\over h}\right)
{\rm d}w\ . \label{Jdef} \end{equation}
Note that $\phi$ is bounded below and $\eta^*$ is positive, and hence 
$J(a,b,\phi,\psi)$ is well-defined.
It then follows from \pref{5upper}, \pref{5aaa} and \pref{Jdef} that
for any dual convex pair of functions $\phi$ and $\psi$, $a\in R^d$ and any $b\in 
\R$,
\begin{equation}  I(F) \ge J(a,b,\phi,\psi)\ .\label{5max} \end{equation}

We let ${\cal U}$ denote the set of all quadruplets $(a,b,\phi,\psi)$ where
$a\in \R^d$,  $b\in \R$, 
and  $\phi$ and $\psi$ are a pair of dual convex functions
with 
\begin{equation} \phi(v) = \infty\qquad{\rm for}\quad |v|>R\ .\label{5cut} \end{equation}
The reason for this last condition is that increasing $\phi$ off of the support 
of $F_0$ can only decrease $\psi$ 
and hence increase
$J$; so we may freely restrict our attention to such dual pairs; see 
\cite{Evans} or \cite{Bren}. This guarantees that
\pref{5lip} holds whenever $(a,b,\phi,\psi)\in {\cal U}$. Indeed, since $\psi$ is 
determined by $\phi$ through the Legendre
transform,
$J$ can be regarded as a functional of $a$, $b$ and $\phi$ alone. However, the 
notation 
with $\phi$ included as a variable is convenient for the exposition.

As we will see below, 
\begin{equation} \min\{I(F) \ | F\in{\cal E}\} = \max\{J( a, b,\phi,\psi)\ |\ ( a, b,\phi,\psi) 
\in {\cal U}\}\ .
\label{MiNmaX} \end{equation}
The parameters $a$ and $b$ will be seen to 
function as Lagrange multipliers guaranteeing that at the maximum on the right, 
$F_1 = \nabla\phi\#F_0$
does belong to~${\cal E}$.

\proclaim{Theorem} There exists $(a_0,b_0,\phi_0,\psi_0)\in{\cal 
U}$ such that
\begin{equation} J(a_0,b_0,\phi_0,\psi_0) \ge J( a, b,\phi,\psi)\label{5win} \end{equation}
for all $(a, b,\phi,\psi)\in{\cal U}$\ .
Furthermore{\rm ,} if 
\begin{equation} F_1(w) = \left(\eta^*\right)'\left({a_0\cdot w + b_0|w|^2/2 + \psi_0(w)\over 
h}\right)\label{dendef} \end{equation}
then $F_1 \in{\cal E}${\rm ,} 
\begin{equation} \nabla \psi_0\# F_1 = F_0\label{906T} \end{equation}
and
\begin{equation} \nabla \psi_0(w) = w + h\nabla\ln(F_1) + hdw - W_2^2(F_0,F_1)\ .\label{thel} \end{equation}
\endproclaim

Note that this gives us a solution of the Euler-Lagrange equation for the 
minimum of $I(F)$ that we derived in the last section. And indeed, since 
$\eta(t) + \eta^*(s) = st$ with 
$$
t= F_1\quad \hbox{and} \quad s = 
{a_0\cdot w + b_0|w|^2/2 + \psi_0(w)\over h}
$$
with $F=F_1$, $\psi = \psi_0$, there is equality in \pref{5young}. 
By \pref{906T}, there is equality in \pref{5upper} when $F=F_1$, $\psi = \psi_0$ 
and $\phi = \phi_0$.
It follows that
$I(F_1) = J(a_0,b_0,\phi_0,\psi_0)$. Together with \pref{5max}, this proves that 
$F_1$ minimizes $I$ 
on ${\cal E}$.
Thus Theorem 5.1 provides us with the minimizer
of the original problem. The advantage of the $J$ functional lies in the 
compactness properties of
the dual convex pairs.  
 
\demo{Proof} First, suppose that the maximizer 
$(a_0,b_0,\phi_0,\psi_0)$ does exist.
Observe that for any real number $\lambda$, 
$(a_0,b_0,\phi_0+\lambda,\psi_0-\lambda)\in{\cal U}$. Then by \pref{5win}
$${{\rm d}\over {\rm 
d}\lambda}J(a_0,b_0,\phi_0+\lambda,\psi_0-\lambda)\bigg|_{\lambda=0} = 0$$
and this clearly leads to
\begin{equation} 1 = \int_{\R^d}\left(\eta^*\right)'
\left({a_0\cdot w + b_0|w|^2/2 + \psi_0(w)\over h}\right){\rm d}w\ .\label{5bnd} \end{equation}
Hence we see that \pref{dendef} does define a probability density.

Next, we shall see below that for some $\varepsilon>0$,
\begin{equation} \int_{\R^d}e^{\varepsilon|w|^2}F_1(w){\rm d}w < \infty\ .\label{5Fer} \end{equation}
This implies that
$$(a,b)\mapsto \int_{\R^d}\left(\eta^*\right)'
\left({a\cdot w + b|w|^2/2 + \psi_0(w)\over h}\right){\rm d}w$$
is a differentiable function of $a$ and $b$ in some neighborhood of $(a_0,b_0)$. 
Assuming this for the moment,
${{\rm d}\over {\rm d}b}J(a_0,b,\phi_0,\psi_0)\big|_{b=b_0} = 0$,
and from this we have that
$${d \over 2} = \int_{\R^d}{|w|^2\over 2}\left(\eta^*\right)'
\left({a_0\cdot w + b_0|w|^2/2 + \psi_0(w)\over h}\right){\rm d}w$$
which means that $F_1$ does indeed satisfy the variance constraint. In the same 
way, differentiating in $a$
shows that $F_1$ does satisfy the mean constraint. Thus, $F_1\in {\cal E}$.

So far, the only variation  made in $\phi_0$, and hence in  $\psi_0$, is 
a shift by an additive constant.
We now let $\zeta$ be any smooth function supported in the interior of $B_R$, 
and define 
$\phi_t = \phi_0+t\zeta$,
and let $\psi_t$ be the Legendre transform of $\phi_t$. While these are not a 
dual pair of convex functions
since $\phi_t$ may fail to be convex, it is nonetheless clear that for all 
sufficiently small $t$,
$J(a_0,b_0,\phi_0,\psi_0) \ge  J(a_0,b_0,\phi_t,\psi_t)$
and thus
$${{\rm d}\over {\rm d}t}J(a_0,b_0,\phi_t,\psi_t)\bigg|_{t=0} = 0\ .$$
As in \cite{Gangbo} $\lim_{t\to 0} (\psi_t(w) - \psi_0(w))/t = 
-\zeta(\nabla\psi_0(w))$ and it follows that
$$
\int_{\R^d}\zeta(v)F_0(v){\rm d}v = \int_{\R^d}\zeta(\nabla\psi_0(w))F_1(w)){\rm 
d}w
,$$
which means that $\nabla\psi_0\#F_1 = F_0\ .$
\vglue4pt

The remaining part of the Euler-Lagrange equation follows from \pref{dendef} by 
simple differentiation:
\begin{equation} h \nabla F_1(w) = \left(a_0 + b_0w + \nabla \psi_0(w)\right)F_1(w)\ 
.\label{5daa} \end{equation}
Hence
$h w\cdot \nabla F_1(w) = \left(a_0\cdot w + b_0|w|^2 + w \cdot 
\nabla\psi_0(w)\right)F_1(w)$,
and integrating both sides we obtain that
$$
b_0 = -(1+h) + {W_2^2(F_0,F_1)\over d}\ .
$$
Even more simply, one sees by integrating \pref{5daa} that $a_0=0$. Thus, 
provided the maximizer exists,
and that $(a,b)\mapsto J(a,b,\phi_0,\psi_0)$ is differentiable in a neighborhood 
of $(a_0,b_0)$, we have
that $F_1\in{\cal E}$, $\nabla\psi_0\#F_1 = F_0$, and  that the Euler-Lagrange 
equation \pref{thel} is satisfied.

To show the existence of an optimizer, we begin by considering any\break
$(a,b,\phi,\psi)$.
We now seek an {\it a priori}  lower bound on $\phi(v)$. Fix any $v\in B_R$ at 
which $\phi$ 
is differentiable.
Then let $w_0 = \nabla\phi(v)$. Since $\psi$ and $\phi$ are dual to one another, 
$v$ belongs to the 
subgradient of $\psi$ at $w_0$, and then by the convexity of $\psi$, for any 
$w\in \R^d$,
$\psi(w) \ge \psi(w_0) + v\cdot(w-w_0)$.

Then since $\psi$ is convex,
and because of the mononicity of $(\eta^*)'$ 
and its specific form, we have that
\begin{eqnarray*}
&&\left(\eta^*\right)'\left({a\cdot w + b|w|^2/2 + \psi(w)\over h}\right) 
\\
&&\qquad \ge \exp\left({(\psi(w_0) - v\cdot w_0)/h}\right)
\left(\eta^*\right)'\left({(a+v)\cdot w + b|w|^2/2\over h}\right)\ .
\end{eqnarray*}
Integrating, and using \pref{5bnd}, we see that $b$ is negative, and obtain
\begin{equation} 1 \ge \exp\left((\psi(w_0)- v\cdot w_0)/h\right)e^{-1}
\exp\left({|a+v|^2 \over 2h|b|}\right)\left({2\pi h\over |b|}\right)^{d/2}\ 
.\hskip.5in\label{5tup} \end{equation}

But $\phi(v) =
-(\psi(w_0)-w_0\cdot v)$ and so
\begin{equation} \phi(v) \ge {|a+v|^2\over 2|b|} - h\left(1 + {d\over 2}\ln\left({|b|\over 2\pi 
h}\right)\right)\ .\label{5tvv} \end{equation}
Integrating against $F_0(v)$, we obtain that
\begin{equation} \int_{\R^d}\phi(v)F_0(v){\rm d}v \ge  {|a|^2\over 2|b|} + {1\over 2|b|} -
h\left(1 + {d\over 2}\ln\left({|b|\over 2\pi h}\right)\right)\ . \hskip.5in\label{5tzz} \end{equation}

Now consider $\tilde\psi$ where 
$\tilde\psi(w) = (1-h)(|w|^2/2) + h\left[1 - (d/2)\ln(2\pi)\right]$, 
so that
$$\int_{\R^d}\left(\eta^*\right)'((-|w|^2/2 + \tilde \psi(w))/h){\rm d}w  = 
\left({1\over 2\pi}\right)^{d/2}\int_{\R^d}e^{-|w|^2/2} =1\ .$$
The dual convex function of $\tilde \psi$ is $\tilde\phi$ where
$$\tilde\phi(w) = 
(|w|^2/2(1-h)) - h\left[1 - (d/2)\ln(2\pi)\right].$$
This does not satisfy \pref{5cut}, and hence
$(0,-1,\tilde\phi,\tilde\psi)$ is not in ${\cal U}$.  However, define 
${\tilde\phi}_R$ by ${\tilde\phi}_R(v)= \tilde \phi(v)$ for $|v|<R$, and 
${\tilde\phi}_R$ by ${\tilde\phi}_R(v)= \infty$ otherwise, and define $\tilde 
\psi_R$ to 
be the dual convex function. Then
$(0,-1,\tilde\phi_R,\tilde\psi_R)$ {\it is} in ${\cal U}$  and
$J(0,-1,\tilde\phi_R,\tilde\psi_R) \ge J(0,-1,\tilde\phi,\tilde\psi)$
since, as we have noted, increasing $\psi$ off the support of $F_0$ 
can only decrease the dual $\psi$, and hence increase $J$. We denote by $J_d(h)$ 
the finite real number  $J(0,-1,\tilde\phi,\tilde\psi),$ depending only on $d$ 
and $h.$ Since it is clear that
$$\sup\{J(a,b,\phi,\psi)\ |\ (a,b,\phi,\psi)\in {\cal U}\} \ge  
J_d(h)\ ,$$
and we seek a maximizer of $J$, we need only consider $(a,b,\phi,\psi) \in 
{\cal U}$ such that
\begin{equation} J(a,b,\phi,\psi) \ge J_d(h)\ .\label{5er} \end{equation}
Furthermore, we may suppose that we have already optimized over
$\phi+\lambda$ and $\psi-\lambda$
so that  \pref{5bnd} holds.  Then from the fact that $(\eta^*)' = \eta^*$, 
$$J(a,b,\phi,\psi) = d-\int_{\R^d}\phi(v)F_0(v){\rm d}v +b {d\over 2} -h\ .$$
In light of this, and \pref{5er}, 
\begin{equation} \int_{\R^d}\phi(v)F_0(v){\rm d}v \le -J_d(h) + d(1+ {b \over 2}) -h \ 
.\label{5uzz} \end{equation}
Combining \pref{5tzz} and \pref{5uzz} we obtain after simplification that 
\begin{equation} 
-J_d(h) + d(1+{b \over 2}) \ge 
{1\over 2|b|} + {a^2\over 2|b|} +{d\over 2}h\ln(2 \pi h) - h{d\over 2}\ln|b|
\ .\hskip.5in\label{55zz} \end{equation}
Recalling that $b$ is negative, it is clear that $|b|$ cannot be too close to 
zero, 
for then the right-hand side becomes greater than $2$. Also, $|b|$ cannot be too 
large, since
as $|b|$ increases, the left-hand side tends linearly to $-\infty$, while the 
right-hand side
only does so logarithmically. Even more  evidently, $|a|$ cannot be too large.

It follows  that there is a constant $c>0$, depending on $h$, so that
\begin{equation} c \le |b| \le 1/c\qquad{\rm and}\qquad |a| < c\ .\label{5abbnd} \end{equation}
Next, use \pref{dendef} to define $F_1$; that is,
\begin{equation} F_1(w) = \left(\eta^*\right)'\left({a\cdot w + b|w|^2/2 + \psi(w)\over 
h}\right)\ .\label{5678} \end{equation}
We may suppose without loss of generality that $a$ and $b$ have been chosen 
optimally so that
$F_1 \in {\cal E}$. Since  $\int_{\R^d}|v|^2F_1(v){\rm d}v =1$,
$$1/2 \le \int_{|w|\le \sqrt{2}}F_1(w){\rm d}w \le 1.\pagebreak$$ 
This together with \pref{5lip} and \pref{5678}  means that for another finite  constant $C$,
\begin{equation} |\psi(w)| \le C + R|w|\label{5lll} \end{equation} 
for all $w$. In particular, with $F_1$ defined as in \pref{5678}, \pref{5Fer}
holds, as claimed.

This gives   all of the {\it a priori} estimates   needed. Consider a 
sequence\hfill\break
$\{(a_n,b_n,\phi_n,\psi_n)\in{\cal U}$, each of which satisfies \pref{5er}. 
First we may optimize in $a_n$ and $b_n$
and carry out the variation over $\phi_n+\lambda$ and $\psi_n-\lambda$.
With these chosen optimally, \pref{5bnd} holds.

Then by the previous paragraphs,
$a_n$ and $b_n$ satisfy \pref{5abbnd} for all $n$.
Passing to a subsequence, we may assume that
$\{a_n\}$ and $\{b_n\}$ converge to the limits $a_0$ and $b_0$ respectively. 

Now for each $n$, define  $F_1^{(n)}$ in terms of $a_n$, $b_n$ and $\psi_n$ 
using \pref{5678}  
Our optimizing sequence is such that  for each $n$, 
$F_1^{(n)}\in {\cal E}$, since, as we have seen, this is what is guaranteed by 
optimality in $a$ and $b$. Moreover, since $a_n$ and $b_n$ satisfy \pref{5abbnd}
for all $n$,
it follows that \pref{5Fer} holds for all $n$ for some fixed $\varepsilon>0$.

Passing to a further subsequence, we have that $\psi_0  = 
\lim_{n\to\infty}\psi_n$ exists uniformly on compact
sets due to \pref{5lll} and the Lipschitz bound. Since for each $n$, $F_1^{(n)}$ 
satisfies \pref{5Fer},
$\lim_{n\to\infty}F_1^{(n)}$ converges strongly in $L^1$. 

It is plain that on $B_R$, passing to a further 
subsequence if need be, we have $\lim_{n\to\infty}\phi_n = \phi$ almost 
everywhere and
$$\lim_{n\to\infty}\int_{B_R}\phi_n(w)F_0(w){\rm d}w =  
\int_{B_R}\phi_0(w)F_0(w){\rm d}w\ .$$ Thus
$J(a_0,b_0,\phi_0,\psi_0) = \lim_{n\to\infty} J(a_n,b_n,\phi_n,\psi_n)$.
Since $\{(a_n,b_n,\phi_n,\psi_n)\}$ was a maximizing sequence, 
$(a_0,b_0,\phi_0,\psi_0)\in {\cal U}$
is the desired maximizer, and all of the properties of $F_1$ and $\psi_0$ 
claimed in the theorem have already 
been shown to be consequences of the corresponding Euler-Lagrange equations. 
\enddemo

Thus, under our given conditions on $F_0$, we have proved the existence of a 
minimizer $F_1$ of $I(F)$.
Now consider an arbitrary element $F_0\in {\cal E}$. Then there exists a convex 
function $\zeta$ on $\R_+$
as in Section 4 such that $\zeta(t)/t$ increases to infinity and
$$\int_{\R^d}\zeta(|v|^2)F_0(v){\rm d}v = C < \infty\ .$$
We approximate $F_0$ in $L^1(\R^d)$ by a sequence of densities $F_0^{(n)}$ such 
that
$$\int_{\R^d}\zeta(|v|^2)F_0^{(n)}(v){\rm d}v < 2C\ $$
for all $n$, and such that for each $n$, $F_0^{(n)}$ is supported in $B_{R_n}$ 
for some radius $R_n$.
Let $F_1^{(n)}$ be the corresponding minimizer of $I(F)$. Then by \pagebreak Theorem 4.2, 
there are numbers $\alpha>0$
and $K<\infty$ so that
\begin{equation} \int_{\R^d}\zeta(\alpha|v|^2)F_1^{(n)}(v){\rm d}v < K\ \label{HKHK} \end{equation}
for all $n$. 

By passing to a subsequence, we may suppose that $F_1^{(n)}$ converges weakly to 
a probability density $F_1$. It is clear that the first moments converge, and by 
\pref{HKHK} it is clear that the second moments converge as well, and hence $F_1 
\in {\cal E}$. Moreover,
since convergence in the $2$-Wasserstein metric is equivalent to weak convergence 
and convergence of the second moments,
$\lim_{n\to\infty}W_2^2(F_1^{(n)},F_1)\break = 0$, and 
$\lim_{n\to\infty}W_2^2(F_0^{(n)},F_0)= 0$. Therefore,
$$\lim_{n\to\infty} W_2^2(F_1,F_0) = W_2^2(F_1^{(n)},F^{(n)}_0)\ .$$
Finally, by weak lower semicontinuity,
$S(F_1) \le \liminf_{n\to\infty} S(F_1^{(n)})$.
It follows that $F_1$ is the minimizer we seek.

Then by dominated convergence, $F_1 = \lim_{n\to\infty}F_1^{(n)}\in{\cal E}$ and 
$F_1$ is the desired minimizer. It is unique by strict convexity.
Thus we have proven the following result:
\proclaim{Theorem} For all $F_0\in {\cal E}${\rm ,} there exists a unique 
$F_1 \in {\cal E}$
such that
$$ I(F_1) \le I(F)$$
for all $F\in {\cal E}${\rm ,} where $I(F)$ is  as defined in {\rm \pref{fudef}.}
\endproclaim
 
We note that on the basis of this result, there is a unique solution to the 
discrete time evolution problem
in which, given initial data $F_0\in {\cal E}$ and a time step $h>0$, $F_n$ is 
defined iteratively in terms of
$F_{n-1}$ by setting $F_n$ to be the minimizer of
$$\left[{W_2^2(  F_{n-1},F)
\over \theta} +hS(F)\right]$$
over ${\cal E}$. We see easily, using the results of Section 4, that if 
we define  
$F^{(h)}(t,v)$ by an appropriate interpolation as in \cite{JKO}, then 
$\lim_{h\to 0}F^{(h)}(t,v) = F(t,v)$
where $F(t,v)$ solves the Fokker-Planck equation 
$${\partial \over \partial t}F(t,v) = \nabla\cdot 
\left(e^{-|v-u|^2/2\theta}\nabla(e^{|v-u|^2/2\theta}F(t,v)\right)
$$
with initial data $F_0$. This equation is of course already well understood, but 
we shall show that this way of
approaching it extends to the nonlinear spatially inhomogeneous kinetic 
Fokker-Planck equation, which is much 
less well understood, \pagebreak
in a related paper.

{\it Open problems}.
We close this section by commenting on two open problems. First, consider the 
variational
problem employed by Jordan, Kinderlehrer and Otto \cite{JKO} to construct 
solutions of the heat equation:
\begin{equation} \inf\{ hS(F) + W_2^2(F, F_0)\}\label{HJHJ}
\end{equation}
in which {\it no constraint} is imposed on the variance of $F$. We conjecture 
that
\begin{equation} \int_{\R^d}|v|^2 F_1(v){\rm d}v > \int_{\R^d}|v|^2 F_0(v){\rm d}v\label{PKPK} \end{equation}
where $F_1$ is the minimizer for \pref{HJHJ}. We can prove this under several 
additional assumptions ---
when $h$ is not too small, when $F_0$ is radial, etc., and we note that if $F_t$ 
solves the heat equation,
\begin{equation} {{\rm d}\over {\rm d}t}\int_{\R^d}|v|^2 F_t(v){\rm d}v = 2d\label{hemon}\end{equation}
for any initial data $F_0$ with finite variance. In Section 3, we have given the 
exact solution of
this variational problem, and we see similar behavior in that case. However, 
we have not been able to prove
\pref{PKPK} in general. It would be most unfortunate if the discrete time problem 
did not possess a
good analog of the basic montonicity property \pref{hemon}, and we do not believe 
that this is the case. 
If \pref{PKPK} were true, it would make it easy to prove Theorem 5.2 by adding on 
a Lagrange multiplier
$\lambda\int_{\R^d}|v|^2 F(v){\rm d}v$
to the functional in \pref{HJHJ}. The existence (and uniqueness) of minimizers 
would follow by the argument
in \cite{JKO} for all $\lambda >0$. Let $F^{(\lambda)}$ denote the minimizer 
corresponding
to a given value of $\lambda \ge 0$. If \pref{PKPK} were true, it would be easy 
to show the existence 
of a value $\lambda_0>0$ for which
$\int_{\R^d}|v|^2 F^{(\lambda_0)}(v){\rm d}v = \int_{\R^d}|v|^2 F_0(v){\rm d}v$.
It would then follow that $F^{(\lambda_0)}$ is the minimizer provided by Theorem 
5.3.

Another open problem concerns the growth of higher moments. 
We note that if $F_t$ solves the heat equation 
for any initial data $F_0$ with zero mean and finite fourth moments,

\centerline{${\displaystyle {{\rm d}\over {\rm d}t}\int_{\R^d}|v|^4 F_t(v){\rm d}v = 12d\theta\ .}$}
\vglue6pt
\noindent
This leads one to hope that if $F_1$ is the minimizer for \pref{HJHJ}, and $F_0$
has zero mean and, say, finite sixth moments, there is a constant 
$C$ depending only on, say, the sixth moments so that
\begin{equation} \int_{\R^d}|v|^4 F_1(v){\rm d}v  \le (1+Ch) \int_{\R^d}|v|^4 F_0(v){\rm d}v\ .\label{PzPK} \end{equation}
This would be helpful in studying the nonlinear kinetic Fokker-Planck 
equation
by these methods. We conjecture that this is true.
We note that to prove \pref{PzPK}, one needs an upper bound on the moments of the 
minimizer $F_1$, while to
prove \pref{PKPK}, one needs a lower bound.
\vglue-12pt

\bye